\documentclass{elsart}
\usepackage[T1]{fontenc}
\usepackage{textcomp,mathrsfs,amssymb,amsmath}
\usepackage[dvips]{graphicx}
\newcommand{\romannum}{\renewcommand{\labelenumi}{\textnormal{(\roman{enumi}).}}}

\usepackage[matrix,arrow]{xy}

\newcommand\vphi{\varphi}
\newcommand\vrho{\varrho}

\newcommand\eps{\varepsilon}
\newcommand\nats{\mathbb{N}}
\newcommand\ints{\mathbb{Z}}
\newcommand\reals{\mathbb{R}}
\newcommand\cplxs{\mathbb{C}}
\newcommand\ball{\mathbb{B}}
\newcommand\sph{\mathbb{S}}
\newcommand\vvoid{\varnothing}
\newcommand\sle{\leqslant}
\newcommand\sge{\geqslant}
\newcommand\limk{\lim\nolimits}
\newcommand\equiva{\Leftrightarrow}
\newcommand\imp{\Rightarrow}
\newcommand\clos[1]{\overline{{#1}}}
\newcommand\var{\llcorner\kern-.3em\lrcorner}

\newcommand\mathtxt[1]{\quad\text{{#1}}\quad}

\newcommand{\nd}{\mathtxt{and}}
\newcommand{\nda}{\ \text{and}\ }

\newcommand\fa{for all }
\newcommand\mathfa[1][{}]{\quad\text{\fa{#1} }}
\newcommand\sth{so that }
\newcommand\scth{such that }

\DeclareMathOperator\dist{dist}

\DeclareMathOperator\ext{ext}
\DeclareMathOperator\id{id}
\DeclareMathOperator\pr{pr}
\DeclareMathOperator\GL{GL}

\makeatletter
\def\Set@Scallop[#1]#2#3{{#1}\Parens{#2}{#3}}
\newcommand\DeclareScalableOperator[2]{%
  \expandafter\def\csname#1\endcsname{\@ifnextchar[{{#2}\Set@Scallop}{{#2}\Set@Scallop[{}]}}
}
\newcommand\Size[7][1]{
                                 \ifx#20%
                                        \def\r@l{}\def\r@m{}\def\r@r{}%
                                 \else%
                                    \ifx#21%
                                           \def\r@l{\bigl}\def\r@r{\bigr}\def\r@m{\bigm}%
                                    \else%
                                           \ifx#22%
                                                 \def\r@l{\Bigl}\def\r@r{\Bigr}\def\r@m{\Bigm}%
                                            \else%
                                                 \ifx#23%
                                                        \def\r@l{\biggl}\def\r@r{\biggr}\def\r@m{\biggm}%
                                                  \else
                                                        \ifx#24%
                                                              \def\r@l{\Biggl}\def\r@r{\Biggr}\def\r@m{\Biggm}%
                                                        \fi%
                                                  \fi%
                                            \fi%
                                      \fi%
                                 \fi%
                                 \ifx#10%
                                       \def\r@m{}%
                                 \fi%
                                 \r@l#3{#4}\r@m#5{#6}\r@r#7%
}%
\newcommand\Set[3]{
                                 \Size{#1}{\{}{#2}{|}{#3}{\}}
}%
\newcommand\Rscp[3]{
                                 \Size[0]{#1}{(}{#2}{|}{#3}{)}%
}%
\newcommand\Parens[2]{
  \Size[0]{#1}{(}{#2}{}{}{)}
}
\newcommand\Bracks[2]{
  \Size[0]{#1}{[}{#2}{}{}{]}
}
\newcommand\Opint[3]{
  \Size[0]{#1}{]}{#2}{,}{#3}{[}%
}
\newcommand\Lopint[3]{
  \Size[0]{#1}{]}{#2}{,}{#3}{]}%
}
\newcommand\Ropint[3]{
  \Size[0]{#1}{[}{#2}{,}{#3}{[}%
}
\newcommand\Norm[2]{
  \Size[0]{#1}{\lVert}{#2}{}{}{\rVert}
}
\newcommand\Abs[2]{
  \Size[0]{#1}{\lvert}{#2}{}{}{\rvert}
}
\newcommand\Span[2]{
  \Size[0]{#1}{\langle}{#2}{}{}{\rangle}
}

\DeclareScalableOperator{Mea}{\mathfrak{M}} 
\DeclareScalableOperator{Ct}{\mathcal{C}} 
\DeclareScalableOperator{Cc}{\mathcal{C}_c} 
\DeclareScalableOperator{Lp}{\mathbf{L}} 
\DeclareScalableOperator{Cst}{\mathrm{C}^*} 
\DeclareScalableOperator{Cred}{\mathrm{C}^*_r} 
\DeclareScalableOperator{Cpt}{\mathbb{K}} 
\DeclareScalableOperator{Lct}{\mathcal{L}} 
\DeclareScalableOperator{closed}{\mathbb F}

\journal{Advances in Mathematics}

\begin{document}

\begin{frontmatter}
  \title{An Index Theorem for Wiener--Hopf Operators}
  \author{Alexander Alldridge\corauthref{aa}}\ead{alldridg@math.upb.de},
  \author{Troels Roussau Johansen\thanksref{trj}}\ead{johansen@math.upb.de}
  \address{Institut f\"ur Mathematik, Universit\"at Paderborn, Warburger Stra\ss{}e, 33098 Paderborn}
  \thanks[trj]{Partly supported by IHP `Harmonic Analysis and Related Problems', HPRN-CT-2001-00273. Currently supported by a DFG post.doc-grant under the International Research \& Training Group ``Geometry and Analysis of Symmetries'', {\texttt{http://irtg.upb.de}}}
  \corauth[aa]{Corresponding author.}  


  \begin{abstract}
    We study the multivariate generalisation of the classical Wiener--Hopf algebra, which is the C$^*$-algebra generated by the Wiener--Hopf operators, given by convolutions restricted to convex cones. By the work of Muhly and Renault, this C$^*$-algebra is known to be isomorphic to the reduced C$^*$-algebra of a certain restricted action groupoid. It admits a composition series, and therefore, a `symbol' calculus. Using groupoid methods, we obtain, in the framework of Kasparov's bivariant $KK$-theory, a topological expression of the index maps associated to these symbol maps in terms of geometric-topological data of the underlying convex cone. This generalises an index theorem by Upmeier concerning Wiener--Hopf operators on symmetric cones. Our result covers a wide class of cones containing the polyhedral and homogeneous cones. 
  \end{abstract}
\end{frontmatter}

\section{Introduction}

  Let $\Omega$ be a closed, pointed, and solid convex cone in the $n$-dimensional real inner product space $X\,$. 
  Wiener--Hopf operators are the bounded operators $W_f$ on $\Lp[^2]0\Omega\,$, defined by
  \[
  W_f\xi(x)=\int_\Omega f(x-y)\xi(y)\,dy\mathfa f\in\Lp[^1]0X\,,\,\xi\in\Lp[^2]0\Omega\,,\,x\in\Omega\ . 
  \]
  Classically, one considers $X=\reals$ and $\Omega=\reals_{\sge0}\,$. Then it is a well-known fact that $1+W_f$ is a Fredholm operator if and only if $1+\Hat f$ is 
  everywhere non-vanishing, where $\Hat f$ denotes Fourier transform.  Then $1+\Hat f$ may be considered as a non-vanishing function on $\sph^1\,$, and its winding number 
  is the index of $1+W_f$ (possibly up to a sign, depending on the normalisation of the Fourier transform).  

  These two facts may be phrased in a modern language, as follows. Let $A$ be the C$^*$-algebra, generated by $W_f\,$, $f\in\Lp[^1]0\reals\,$. Then the map 
  $W_f\mapsto\Hat f$ extends to a morphism of C$^*$-algebras, and we have a short exact sequence
  \[
  \xymatrix{%
    0\ar[r]&\mathbb K\Parens0{\Lp[^2]0{0,\infty}}\ar[r]&A\ar[r]&\Ct[_0]0\reals\ar[r]&0\ ,}
  \]
  $\mathbb K$ denoting the set of all compact operators and $\mathcal C_0$ that of the continuous functions vanishing at infinity. Moreover, the index of operators with 
  given (matrix) symbols corresponds to the connecting map in operator $K$-theory, $K^1_c(\reals)=K_1(\Ct[_0]0\reals)\to K_0(\mathbb K)=\ints\,$, induced by this short 
  exact sequence. Computing the index of Fredholm Wiener--Hopf operators (in fact, Fredholm matrices in the unitisation of $A$) amounts to the computation of this map. 

  The generalisation of these results to multivariate situations, where $X$ has dimension $n\sge2$ and $\Omega$ is a closed convex cone, has been a long-standing problem, 
  and various authors have contributed to it over time. We shall only give the very briefest overview at this point and refer to the induction of \cite{alldridge-johansen-wh1} 
  and the literature cited there for a more in-depth account of these historical matters. 

  In the multivariate case, the above short exact sequence is replaced by a composition series, i.e.~a filtration by closed ideals whose subquotients are continuous trace 
  algebras, and the corresponding index maps are the $d^1$ differentials of the Atiyah--Hirzebruch homology spectral sequence in $K$-theory induced by this filtration. A 
  composition series was obtained in the opposite case of polyhedral and symmetric (i.e.~homogeneous and self-dual) cones by Muhly and Renault \cite{muhly-renault-wh}; at this 
  point, no index theorems were established. Upmeier \cite{upmeier-wh} gave an alternative description of the composition series in the symmetric case and computed 
  the index maps. (In fact, he considers the C$^*$-algebras of Toeplitz operators on bounded symmetric domains and develops the corresponding theory for them. If the 
  bounded symmetric domain has a Siegel realisation as a tube domain, then the Wiener--Hopf algebra $A$ occurs as an ideal in the Toeplitz algebra, as in the one-variable 
  case.) 

  This makes the class of symmetric cones the only class of cones for which a complete index theory of multivariate Wiener--Hopf operators has been developed. One obstacle 
  to the generalisation of the finer $K$-theoretical results beyond symmetric cones seemed to be that Upmeier's work relies 
  on the use of the theory of Jordan algebras, which is only available in this case. Due to the refinement of the theory of groupoid C$^*$-algebras and operator $K$ theory 
  over the last 20 years, we have been able to follow the groupoid approach of Muhly and Renault, thereby obtaining an independent proof of Upmeier's results and 
  simultaneously extending them to a large class of cones which, in particular, allows the treatment of polyhedral and homogeneous (not necessarily symmetric) cones by 
  the same methods. 

  Let us explain this approach. The C$^*$-algebra generated by the Wiener--Hopf operators $W_f\,$, $f\in\Lp[^1]0X\,$, is isomorphic to the reduced 
  groupoid C$^*$-algebra $\Cred0{\mathcal W_\Omega}$ of the `Wiener--Hopf groupoid' $\mathcal W_\Omega\,$.

  The groupoid $\mathcal W_\Omega$ is a certain restriction of an action groupoid, constructed as follows. The vector space $X$ acts by translation on its space of closed 
  subsets $\closed0X\,$, which is a compact metrisable space when endowed with the Fell topology. The set $-\Omega$ is a point of $\closed0X\,$, and the closure $\clos X$ 
  of the orbit through this point under the action of $X$ is compact. The orbit closure for the action of the subsemigroup $\Omega\subset X$ is denoted by 
  $\clos\Omega\,$. Now, if $\clos X\rtimes X$ denotes the usual action groupoid given by the action of the group $X$ on the space $\clos X\,$, the Wiener--Hopf groupoid 
  is defined as the restriction $\mathcal W_\Omega=(\clos X\rtimes X)|\clos\Omega$ to the non-invariant subset $\clos\Omega\,$. We refer to 
  \cite{hn-wiener1,alldridge-johansen-wh1} for the details of the construction.

  The groupoid $\mathcal W_\Omega$ describes a lamination of the space $\clos\Omega\,$, and in \cite{alldridge-johansen-wh1}, we have constructed 
  a composition series of the algebra $\Cred0{\mathcal W_\Omega}$ associated to the orbit type stratification for the groupoid's action, which we proceed to describe. 

  Recall that a subset $F$ of a convex set $C$ is a \emph{face} if any segment in $C$ whose midpoint belongs to $F$ lies completely in $F\,$. Order the 
  dimensions of the faces of the dual cone $\Omega^*$ increasingly, 
  \[
  \{0=n_0<n_1<\dotsm<n_d=n\}=\Set1{\dim F}{F\subset\Omega^*\text{ face}}\ .
  \]
  Let $P_j$ be the set of $n_{d-j}$-dimensional faces of $\Omega^*\,$. Any element of $P_j\,$, considered as a point of $\clos\Omega\,$, has the same orbit type. Here, recall 
  that an \emph{orbit type} of a groupoid $\mathcal G$ with structure maps $r$ and $s$ is the conjugacy class, under the action of $\mathcal G\,$, of an isotropy group 
  $\mathcal G(x)=r^{-1}(x)\cap s^{-1}(x)\,$, $x\in\mathcal G^{(0)}\,$, of $\mathcal G\,$. The isotropies of the groupoid $\mathcal W_\Omega$ are linear subspaces of $X$ and
  their conjugacy is classified by dimension. The orbit type of every point of $P_j$ then turns out to be that of an $(n-n_{d-j})$-dimensional linear subspace of $X\,$. 
  
  Now, let $Y_j\subset\clos\Omega$ be the `orbit type stratum' of $P_j\,$, i.e.~the set of all points of $\clos\Omega$ with the same orbit type as the points of $P_j\,$. 
  Moreover, \emph{assume that the sets $P_j$ are compact} \fa $j\,$, in the space of all closed subsets of $X\,$,  
  endowed with the Fell topology. Then the orbit type strata $Y_j$ cover $\clos\Omega\,$, and naturally form fibre bundles over the spaces $P_j\,$. 

  We call cones whose sets of faces of fixed dimension are compact \emph{facially compact}. This condition obtains for $\Omega^*$ if $\Omega$ is, e.g., polyhedral or 
  symmetric (i.e.~homogeneous and self-dual). In these special cases, the sets $P_j$ are, respectively, finite sets and certain compact homogeneous spaces including, 
  in particular, all spheres.

  The strata $Y_j$ are closed and invariant, and the unions $U_j=\bigcup_{i=0}^{j-1}Y_i$ are open and invariant. Thus, we obtain ideals
  $I_j=\Cred0{\mathcal W_{\Omega}|U_j}$ of the Wiener--Hopf C$^*$-algebra $\Cred0{\mathcal W_\Omega}\,$, and extensions
  \[
  \xymatrix@C-16pt{%
    0\ar[r]&\Cred0{\mathcal W_\Omega|Y_j}\ar[r]&I_{j+1}/I_{j-1}=\Cred0{\mathcal W_\Omega|(U_{j+1}\setminus U_{j-1})}\ar[r]&\Cred0{\mathcal W_\Omega|Y_j}\ar[r]&0}\ .
  \]
  Moreover, the fibre bundle structure of the strata induces Morita equivalences $\mathcal W_\Omega|Y_j\sim\Sigma_j$ where $\Sigma_j$ is 
  the `co-tautological' vector bundle over the space $P_j$ whose fibre at the face $F$ is the orthogonal complement $F^\perp\,$. (Alternatively, $\Sigma_j$ is the 
  restriction of the isotropy group bundle of $\mathcal W_\Omega$ to the section $P_j$ of the bundle $Y_j\to P_j\,$.) Thus, $\Cred0{\mathcal W_\Omega}$ is solvable of 
  length $d\,$, and its spectrum can be computed in terms of a suitable gluing of the bundles $\Sigma_j\,$. As a particular case, one obtains the classical Wiener--Hopf 
  extension (associated to $X=\reals$ and $\Omega=\reals_{\sge0}$).

  The above extensions induce index maps $K^1_c(\Sigma_j)\to K^0_c(\Sigma_{j-1})\,$, given as the Kasparov product with the $KK$ class
  \[
  \partial_j\in KK^1(\Cred0{\mathcal W_\Omega|Y_j},\Cred0{\mathcal W_\Omega|Y_{j-1}})
  \]
  representing the above extension. As already noted, these maps may also be viewed as the $d^1$ 
  differentials of the Atiyah--Hirzebruch type homology spectral sequence induced by the filtration. In \cite{alldridge-johansen-wh1}, we expressed 
  the index map $\partial_j$ as the family index of certain continuous Fredholm families of operators on a continuous field of Hilbert spaces over $\Sigma_{j-1}\,$.

  In this paper we prove a formula for $\partial_j$ which expresses the latter through topological data. This formula generalises Upmeier's 
  \cite{upmeier-wh} result for symmetric cones to a larger class of convex cones. Let us stress that our methods are different from his, and originate in the 
  non-commutative geometry approch to foliation C$^*$-algebras, benefiting from the recent progress in this field. Our result applies both to cones with a 
  large number of automorphisms (such as homogeneous cones which need not be symmetric), and those with only few (such as polyhedral cones). 
  
  We proceed to describe our index formula. To stress the analogy, we deliberately adopt Upmeier's notation, despite the differences in our assumptions and methods.
  Assume that the cone $\Omega$ has a facially compact and locally smooth dual cone (compare 
  section 6 for the definition of `locally smooth').  Consider the compact space $\mathcal P_j$ consisting of all pairs $(E,F)\in P_{j-1}\times P_j$ \scth $E\supset F\,$. 
  It has projections
  \[
  \xymatrix@C+8ex{P_{j-1}&\mathcal P_j\ar[l]_-{\xi}\ar[r]^-{\eta}&P_j}
  \]
  which need not be surjective unless $j=1,d$ (although they are in the polyhedral and symmetric cases). The projection $\xi:\mathcal P_j\to P_{j-1}$ turns
  $\mathcal P_j$ into a fibrewise $\mathcal C^1$ manifold over the compact base $\xi(\mathcal P_j)\,$. Moreover, $\eta^*\Sigma_j$ is the trivial line bundle over
  $T\mathcal P_j\oplus\xi^*\Sigma_{j-1}$ if $T\mathcal P_j$ denotes the fibrewise tangent bundle. Then we have the following theorem

  \begin{thm}
    The $KK$-theory element $\partial_j$ representing the $j^{\text{th}}$ Wiener--Hopf extension is given by
    \[
    \pmb\zeta_*\partial_j=\pmb\eta^*\Bracks0{y\otimes\tau_j}
    \]
    where $y\in KK^1(\cplxs,S)$ represents the classical Wiener--Hopf extension, $\pmb\eta$ is associated to the projection $\eta^*\Sigma_j\to\Sigma_j\,$, and $\pmb\zeta$ is
    associated to the inclusion $\Sigma_{j-1}|\xi(\mathcal P_j)\subset\Sigma_{j-1}\,$. Here, 
    \[
    \tau_j\in KK(\Cred0{T\mathcal P_j\oplus\xi^*\Sigma_{j-1}},\Cred0{\Sigma_{j-1}|\xi(\mathcal P_j)})
    \]
    represents the Atiyah--Singer family index for $T\mathcal P_j\oplus\xi^*\Sigma_{j-1}\,$, considered as a vector bundle over $\Sigma_{j-1}|\xi(\mathcal P_j)\,$.
  \end{thm}

  To illustrate, we first consider the special case $j=d\,$. Here, $\eta$ is constant ($P_d=\{0\}\,$, $\Sigma_d=X$), $\xi$ is the identity, and in particular,
  surjective. The fibres of $\xi$ are points, so $T\mathcal P_{d-1}=0\,$. The vector space $X$ is turned into the trivial real line bundle over $\Sigma_{d-1}$ by
  letting the fibre at $(E,u)\in\Sigma_{d-1}$ be the line spanned by the extreme ray $E$ of $\Omega^*\,$. We have that $\tau_{d-1}$ is the identity, so our index
  formula in this case is just $\partial_d=\pmb\eta^*y\,$, which recovers the case of the classical Wiener--Hopf extension for $\Omega=\reals_{\sge0}\,$. 

  A more interesting special case is $j=1\,$. Here, $P_0=\{\Omega^*\}$ is the point, and $\Sigma_0$ the zero bundle over the point. So, $\xi$ is constant, and $\eta$ is
  the identity. The set $\mathcal P_1=P_1$ consists of all maximal-dimensional proper faces of $\Omega^*\,$. Their dual faces $\check F$ are exposed extreme rays of
  $\Omega\,$. The tangent bundle $TP_1$ has at the face $F$ the fibre $F^\perp\cap\check F^\perp\,$. It is important to note that for non-polyhedral cones, this space is
  usually non-zero, the simplest case being that of the three-dimensional Lorentz cone, where $TP_1$ is the tangent bundle of $\sph^1\,$. In any case, for $j=1\,$, 
  $\tau_1\in KK(TP_1,\cplxs)=K_0(TP_1)$ induces on $K$-theory the Atiyah--Singer index $K^0_c(TP_1)\to\ints$ for 
  $TP_1\,$, and we have $\partial_1=\pmb\eta^*(y\otimes\tau_1)\,$.

  In a forthcoming paper, the first-named author shows that for a polyhedral cone $\Omega\,$, the complex formed by the index maps $\partial_j$ and the 
  cohomology groups of the vector bundles $\Sigma_j$ is precisely the augmented cellular complex associated to the CW complex given by a polygonal section 
  of $\Omega\,$. This suggests an interesting connection between the combinatorics of polyhedra and the $K$-theory of subquotients of the Wiener--Hopf algebra. 
  
  Let us sketch our strategy of proof. We first observe that $KK^1$ elements representing extensions induced by restrictions of groupoids to open invariant subsets 
  behave naturally in the category $KK$ under proper groupoid homorphisms, even if these do not induce $*$-morphisms of the groupoid C$^*$-algebras. 

  We then construct an appropriate proper homomorphism $\vphi$ to the restricted Wiener--Hopf groupoid $\mathcal W_\Omega|(U_{j+1}\setminus U_j)$ which defines the 
  extensions giving rise to the index maps $\partial_j\,$. The domain of $\vphi$ is the fibred product of the classical Wiener--Hopf groupoid 
  $\mathcal W_{\reals_{\sge0}}$ and the tangent groupoid of a fibrewise $\mathcal C^1$ groupoid $\mathcal D_j\,$. Thus, the idea is that the extensions $\partial_j$ are 
  classical Wiener--Hopf extensions with a (possibly highly non-trivial) `twist'. 

  Using naturality, we pull back $\partial_j$ along $\vphi\,$. As is to be expected, the result is a Kasparov product $y\otimes\tau$ where $y$ represents the classical 
  Wiener--Hopf extension and $\tau=\tau_j$ is the `Connes--Skandalis map' associated to the tangent groupoid. It remains to express $y\otimes\tau$ by topological 
  means, but this can be done essentially by standard procedures. 
  
  To summarise, the construction of the groupoid $\mathcal D_j$ and the homomorphism $\vphi$ constitutes the main step of the proof; the remainder of our paper 
  consists in building up a little machinery.

  Following the basic philosophy that an index formula should be the consequence of a naturality of transformations, applied to a 
  particular homomorphism, we have organised our material as follows. Sections 2 to 5 consist in collecting tools which are applied only in section 6 to prove our main 
  theorem. Here, sections 2 and 4 contain to our knowledge new constructions, whereas sections 3 and 5 contain extensions of known techniques. 

  More precisely, in Section 2, we treat the naturality of extensions by expressing the mapping cone construction for groupoid C$^*$-algebras by a groupoid construction. 
  Section 3 concerns fibrewise differentiable groupoids. After recalling basic definitions, we study the tangent groupoid and introduce the (fibrewise) Connes--Skandalis 
  map $\tau\,$. In section 4, we construct the `cone' $\mathbb W\mathcal G$ over the tangent groupoid of a given fibrewise $\mathcal C^1$ groupoid $\mathcal G\,$, and use 
  it to compute $y\otimes\tau$ as the extension of a groupoid C$^*$-algebra; here, the naturality of extensions (from section 2) also enters. We effect the computation 
  of $\tau$ in topological terms in section 5 by adapting Connes's familiar construction of the classifying space in  for the tangent groupoid of a manifold. In Section 6, 
  we finally consider the Wiener--Hopf groupoid. We construct the fibrewise $\mathcal C^1$ groupoid $\mathcal D_j\,$, and a proper 
  homomorphism $\mathbb W\mathcal D_j\to\mathcal W_\Omega|(U_{j+1}\setminus U_{j-1})\,$. An application of the methods previously established proves the index formula.

\section{Groupoid Extensions and Naturality}

\subsection{Preliminaries}\label{preliminaries-section}

  We will consider groupoids $\mathcal G\,$, usually locally compact (and Hausdorff), and often endowed with a (continuous left) Haar system 
  $(\lambda^x)_{x\in\mathcal G^{(0)}}\,$. In the latter case, we will consider the full and reduced groupoid C$^*$-algebras $\Cst0{\mathcal G}$ and $\Cred0{\mathcal G}\,$, 
  where we suppress the Haar system from the notation. We will use these concepts freely, and refer the reader to \cite{renault-groupoid-cst} for details.

  To fix our notation and terminology, we collect some well-known facts on generalised morphisms and related matters. In what follows, let $\mathcal G,\mathcal H$ 
  be locally compact (Hausdorff) groupoids whose source and range maps are open.

  A \emph{generalised morphism} from $\mathcal G$ to $\mathcal H$ is a locally compact space $Z\,$, together with maps
  $\xymatrix@1{\mathcal G^{(0)}&Z\ar[l]_-{r}\ar[r]^-{s}&\mathcal H^{(0)}}$ \scth $\mathcal G$ acts from the left on $Z$ relative $r\,$, $\mathcal H$
  acts from the right on $Z$ relative $s\,$, the actions commute, and $r:Z\to\mathcal G^{(0)}$ is a principal $\mathcal H$ fibration (which is to say that $\mathcal H$ acts
  properly and freely on $Z\,$, transitively on the fibres of $r\,$, and $r$ is open and surjective). To fix terminology, a continuous homomorphism (i.e., a functor)
  $\mathcal G\to\mathcal H$ will be called a \emph{strict morphism}. For these definitions, compare
  \cite{tu-xu-laurentgengoux-twistedkthy,tu-groupoidcohom,legall-equivarkasparov,hilsum-skandalis,connes-skandalis-longindex}.

  Strict morphisms $f:\mathcal G\to\mathcal H$ induce generalised morphisms $Z_f\,$, defined by $Z_f=\mathcal G^{(0)}\times_{\mathcal H^{(0)}}\mathcal H\,$, 
  with respect to $f:\mathcal G^{(0)}\to\mathcal H^{(0)}\,$. The action of $\mathcal H$ is the obvious one, and the action of $\mathcal G$ is given by
  \[
  \gamma.(x,\eta)=\Parens0{r(\gamma),f(\gamma)\eta}\mathfa\gamma\in\mathcal G^x\,,\,\eta\in\mathcal H_{f(x)}\ .
  \]
  Composition of generalised morphisms $Z:\mathcal G\to\mathcal H\,$, $W:\mathcal H\to\mathcal I$ is given by
  \[
  W\circ Z=Z\times_{\mathcal H}W=\Parens0{Z\times_{\mathcal H^{(0)}}W}/\mathcal H\ ,
  \]
  where the action of $\mathcal H$ is diagonal: 
  \[
  (z,w).\tau=(z\tau,\tau^{-1}w)\mathtxt{whenever}s(z)=r(\tau)=r(w)\ .
  \]
  This composition is compatible with the composition of strict morphisms, up to canonical isomorphism. Locally compact groupoids with equivalence classes of generalised 
  morphisms form a category (without passing to isomorphism classes, one obtains a $2$-category); its isomorphisms are called \emph{Morita equivalences}. In the sequel, we 
  shall be somewhat lax in our use of terminology, and will not distinguish between generalised morphisms and their equivalence classes.

  Following \cite{tu-nonhausgrp}, we shall say that a generalised morphism $Z$ is \emph{locally proper} if the action of $\mathcal G$ is proper, and \emph{proper} if in 
  addition, all inverse images of compacts under $s:Z\to\mathcal H^{(0)}$ are $\mathcal G$-compacts. Equivalently, the induced map
  $\mathcal G\backslash Z\to\mathcal H^{(0)}$ is proper. The composition of proper generalised morphisms is proper, and any Morita equivalence is proper.

  Generalising the construction of an equivalence bimodule given by Muhly--Renault--Williams \cite{muhly-renault-williams-equivalences}, Tu
  \cite[Theorem 7.8, Remark 7.17]{tu-nonhausgrp} has associated to any proper generalised morphism $Z:\mathcal G\to\mathcal H$ of locally compact groupoids with Haar
  systems a trivially graded right $\Cred0{\mathcal H}$-Hilbert module, on which $\Cred0{\mathcal G}$ acts by compact endomorphisms. In particular, this defines 
  a Kasparov cycle $KK(Z)\in KK(\Cred0{\mathcal H},\Cred0{\mathcal G})\,$.

  This correspondence is cofunctorial in the following sense. Let $KK$ be the category whose objects are separable C$^*$-algebras, and whose $\mathrm{Hom}$ functor is 
  $KK(-,-)\,$, with composition given by the Kasparov product. The map which associates to each locally compact groupoid with Haar system its reduced groupoid
  C$^*$-algebra, and to each proper generalised morphism $Z$ as above the cycle $KK(Z)\,$, is a cofunctor. In particular, Morita equivalences give rise to $KK$ equivalences.

\subsection{The Mapping Cone Groupoid}\label{short-ex-seq}
  Let $\mathcal G$ be a locally compact $\sigma$-compact (Hausdorff) groupoid with Haar system $(\lambda^u)$ and $U\subset\mathcal G^{(0)}$ an open invariant subset. Set
  $F=\mathcal G^{(0)}\setminus U\,$. As is well-known \cite[Proposition 4.5]{renault-groupoid-cst}, \cite[2.4]{hilsum-skandalis}, \cite[Proposition 2.4.2]{ramazan-thesis}, there 
  is a short exact sequence
  \begin{equation}
    \label{eq:grp-ext}\xymatrix{0\ar[r]&{\Cst0{\mathcal G|U}}\ar[r]^-j&{\Cst0{\mathcal G}}\ar[r]^-q&{\Cst0{\mathcal G|F}}\ar[r]&0}
  \end{equation}
  where $j$ is given by extension of compactly supported functions by zero, and $q$ is the integrated version of the proper homomorphism given by the inclusion
  $\mathcal G|F\subset\mathcal G\,$. In particular, the $KK$ theory class of $q$ is the $KK$ theory class induced, via the functor $KK(-)$ introduced above, by this inclusion. 

  Moreover, if the groupoid $\mathcal G|F$ is topologically amenable, the corresponding sequence of reduced groupoid C$^*$-algebras
  \[
  \xymatrix{0\ar[r]&{\Cred0{\mathcal G|U}}\ar[r]^-j&{\Cred0{\mathcal G}}\ar[r]^-q&{\Cred0{\mathcal G|F}}\ar[r]&0}
  \]
  is also exact. 

  As does any extension of C$^*$-algebras, these sequences induce certain $KK^1$ elements, commonly constructed by considering mapping cones. 
  To better understand these, we will in this section describe the mapping cones for the above sequences as the full or reduced C$^*$-algebra of a groupoid.
  
  For brevity, we denote $I(U)=\Cst0{\mathcal G|U}\,$, $I_r(U)=\Cred0{\mathcal G|U}\,$. Within the corresponding groupoid C$^*$-algebras of $\mathcal G\,$,
  these are given as the closures of the image of $\Cc0U\,$. We briefly recall that $I(-)$ and, under suitable conditions, $I_r(-)\,$, respect intersections and unions. This is 
  probably well-known; for the lack of a reference and to the reader's convenience, we give a proof.

\begin{prop}\label{ideals-opinv-lattice}
  Let $U_1,U_2\subset\mathcal G^{(0)}$ be open and invariant. Then
  \[
  I(U_1)+I(U_2)=I(U_1\cup U_2)\nd I(U_1)\cap I(U_2)=I(U_1\cap U_2)\ .
  \]
  Moreover, $I_r(U_1)+I_r(U_2)=I_r(U_1\cup U_2)\,$, and if $\mathcal G|(U_1\setminus U_2)$ is topologically amenable, then $I_r(U_1)\cap I_r(U_2)=I_r(U_1\cap U_2)\,$.
\end{prop}

\begin{pf}
  Since $\mathcal G^{(0)}$ is normal, by Urysohn's lemma, there exists a partition of unity subordinate to the cover $(U_1,U_2)\,$. Thus, 
  $\Cc0{U_1}+\Cc0{U_2}=\Cc0{U_1\cup U_2}\,$, and there are canonical isomorphisms
  \begin{align*}
  I(U_1)/I(U_1)\cap I(U_2)&\cong\Parens0{I(U_1)+I(U_2)}/I(U_2)\\
  &=I(U_1\cup U_2)/I(U_2)\cong\Cst0{\mathcal G|(U_1\setminus U_2)}\ .
  \end{align*}
  Thus, we have a commutative diagram with exact lines,
  \[
  \xymatrix{%
    0\ar[r]&I(U_1)\cap I(U_2)\ar[r]&I(U_1)\ar[r]\ar@{=}[d]&{I(U_1)/\Parens0{I(U_1)\cap I(U_2)}}\ar[r]\ar[d]&0\\
    0\ar[r]&I(U_1\cap U_2)\ar[r]&I(U_1)\ar[r]&{\Cst0{\mathcal G|U_1\setminus U_2}}\ar[r]&0}
  \]
  Since the vertical arrows are isomorphisms, the kernels of the rightmost nonzero horizontal arrows coincide, so $I(U_1)\cap I(U_2)=I(U_1\cap U_2)\,$. If the groupoid 
  $\mathcal G|(U_1\setminus U_2)$ is topologically amenable, the same argument applies on the level of reduced groupoid C$^*$-algebras.
\end{pf}

  We wish to express the $KK^1$ element associated to an extension of groupoid C$^*$-algebras in groupoid terms. To that end, we recall the usual construction of the 
  first Puppe sequence in $KK$ theory. Given a $*$-morphism $q:A\to A''$ of separable C$^*$-algebras, the mapping cone $C_q$ is defined as the pullback of
  \[
  \xymatrix{A\ar[r]^-{q}&A''&CA''\ar[l]_-{e_0}}
  \]
  where $CA''=\Ct[_0]0{\Ropint001}\otimes A''$ and $e_0(f)=f(0)\,$. The diagram
  \[
  \xymatrix{A&SA''\ar[l]_-{0}\ar[r]^-{\subset}&CA''}
  \]
  induces a map $SA''\to C_q\,$, where we write $SA=\Ct[_0]0{\Opint001}\otimes A$ for the suspension; and the sequence
  \[
  \xymatrix{SA''\ar[r]&C_q\ar[r]&A\ar[r]^-{q}&A''}
  \]
  is called a \emph{mapping cone triangle}. By applying $KK(B,-)\,$, it gives rise to a long exact sequence of Abelian groups\ ,
  \[
  \xymatrix@C+2ex@R-3ex{%
    &\ar@{-->}[r]&KK(B,SA)\ar[r]^-{KK(B,Sq)}&KK(B,SA'')\ar`[r]`[dlll]`[ddlll]`[ddll][ddll]&\\
    &&&&\\
    &KK(B,C_q)\ar[r]&KK(B,A)\ar[r]^-{KK(B,q)}&KK(B,A'')\ar@{-->}[r]&}
  \]
  called the first Puppe sequence \cite[Theorem 19.4.3]{blackadar-ktheory}.  

  If $q$ is the quotient map of a semi-split extension (e.g.~if $A''=A/A'$ is nuclear), then $C_q$ and $A'=\ker q$ are $KK$-equivalent \cite[Theorem 19.5.5]{blackadar-ktheory} via
  \[
  \text{the map}\quad A'\to C_q\quad\text{induced by}\quad\xymatrix{A&A'\ar[l]_-{q}\ar[r]^-{0}&CA''}\ .
  \]
  Thus, up to a $KK$-equivalence, the connecting map $\partial:KK(B,SA'')\to KK(B,A')$ is given by application of the
  functor $KK(B,-)$ to the natural map $SA''\to C_q\,$, cf.~\cite[Theorem 19.5.7]{blackadar-ktheory}. (In fact, more precisely, to its $KK$-theory class in $KK(SA'',C_q)\,$, but
  this amounts to the same \cite[Proposition 18.7.1]{blackadar-ktheory}.) It is also given by the Kasparov product with the element representing the extension.

  From this discussion, it appears to be desirable to give a description of $C_q$ as a groupoid C$^*$-algebra in the case of the extension (\ref{eq:grp-ext}). To that end, consider 
  the embedding $(0,\id):\mathcal G|F\to\Ropint001\times\mathcal G|F$ and form $\mathcal C_F=\mathcal G\cup_{\mathcal G|F}\Parens0{\Ropint001\times\mathcal G|F}\,$.
  Since $\mathcal G$ is the complement in $\mathcal C_F$ (of the open subset $\Opint001\times\mathcal G|F\,$, and $\Ropint001\times\mathcal G|F$ is the complement of the open subset
  $\mathcal G|U\,$, $\mathcal G$ and $\Ropint001\times\mathcal G|F$ are closed in $\mathcal C_F\,$.

  Let $\mathcal C_F^{(0)}=\mathcal G^{(0)}\cup_F\Parens0{\Ropint001\times F}\,$, and let $r,s:\mathcal C_F\to\mathcal C_F^{(0)}$ be induced by
  \[
  \xymatrix@C=30pt{%
    \mathcal G\ar@<0.5ex>[r]^-{r}\ar@<-0.5ex>[r]_-{s}&\mathcal G^{(0)}\ar[r]&\mathcal C_F^{(0)}%
    &{\Ropint001\times F}\ar[l]&{\Ropint001\times\mathcal G|F}\ar@<-0.5ex>[l]_-{r}\ar@<+0.5ex>[l]^-{s}}\ .
  \]
  Then there is a continuous bijection $\mathcal G^{(2)}\cup_{(\mathcal G|F)^{(2)}}\Parens0{\Ropint001\times\mathcal G}^{(2)}\to\mathcal C_F^{(2)}$ where as usual 
  $\mathcal G^{(2)}=\mathcal G\times_{\mathcal G^{(0)}}\mathcal G\,$, etc. The images of $\mathcal G^{(2)}$ and $\Parens0{\Ropint001\times\mathcal G|F}^{(2)}$ 
  in $\mathcal C_F^{(2)}$ are closed, \sth this map is in fact a homeomorphism. By this token, the operations of $\mathcal G$ and $\Ropint001\times\mathcal G|F$ 
  induce operations on $\mathcal C_F\,$, making the latter a topological groupoid.

  Finally, the inclusions $\mathcal G\subset\mathcal C_F\supset\Ropint001\times\mathcal G|F$ being proper, we have continuous maps
  \[
  \xymatrix@C=30pt{\mathcal G^{(0)}\ar[r]^-{(\lambda^u)}&{\Mea0{\mathcal G}}\ar[r]%
    &{\Mea0{\mathcal C_F}}&{\Mea0{\Ropint001\times\mathcal G}}\ar[l]&{\Ropint001\times F}\ar[l]_-{(\delta_t\otimes\lambda^u)}}
  \]
  of the spaces of Radon measures, endowed with their $\sigma(\mathfrak M,\mathcal C_c)$-topologies. Since these maps coincide on $F\,$, they induce a continuous
  map $\mu:\mathcal C_F^{(0)}\to\Mea0{\mathcal C_F}$ which can be seen to define a Haar system. We have established the following proposition.

\begin{prop}
  The space $\mathcal C_F$ is a locally compact $\sigma$-compact groupoid with a Haar system, and the inclusions
  $\mathcal G\subset\mathcal C_F\supset\Ropint001\times\mathcal G|F$ are proper homomorphisms.
\end{prop}

\begin{thm}
  For the quotient map $q$ in the extension (\ref{eq:grp-ext}), and its mapping cone $C_q\,$, we have $C_q\cong\Cst0{\mathcal C_F}\,$. If $\mathcal G|F$ is 
  topologically amenable, and $q$ is the quotient map of the corresponding short exact sequence of reduced groupoid C$^*$-algebras, then the associated mapping cone 
  is $C_q\cong\Cred0{\mathcal C_F}\,$.
\end{thm}

\begin{pf}
  Note that $\Cst0{\Ropint001\times\mathcal G|F}\cong\Ct[_0]0{\Ropint001}\otimes\Cst0{\mathcal G|F}\,$. Hence, the commutative square of proper homomorphisms
  \[
  \xymatrix{%
    \mathcal G|F\ar[r]\ar[d]&{\Ropint001\times\mathcal G|F}\ar[d]\\
    \mathcal G\ar[r]&\mathcal C_F}
  \]
  preserves Haar systems, and thus integrates to the commutative square of $*$-morphisms
  \[
  \xymatrix{%
    {\Cst0{\mathcal G|F}}&{C\Cst0{\mathcal G|F}}\ar[l]_-{e_0}\\
    {\Cst0{\mathcal G}}\ar[u]^-{q}&{\Cst0{\mathcal C_F}}\ar[u]_-{p_2}\ar[l]^-{p_1}}
  \]
  Thus, we obtain a $*$-morphism $p:\Cst0{\mathcal C_F}\to C_q$ whose kernel is $\ker p_1\cap\ker p_2\,$.

  Now, $\mathcal G=\mathcal C_F|\mathcal G^{(0)}$ and $\Ropint001\times\mathcal G|F=\mathcal C_F|\Parens0{\Opint001\times F}\,$, so $p_1$ resp.~$p_2$ are
  the quotient maps for extensions of type (\ref{eq:grp-ext}) for the groupoid $\mathcal C_F$ and the open invariant subsets 
  $\mathcal C_F^{(0)}\setminus\mathcal G^{(0)}\!=\,\Opint001\times F$ and $\mathcal C_F^{(0)}\setminus\Parens0{\Ropint001\times F}=U\,$, respectively. Thus, 
  \[
  \ker p_1\cap\ker p_2=I(\Opint001\times F)\cap I(U)=I(\vvoid)=0\ ,
  \]
  by Proposition \ref{ideals-opinv-lattice}.

  As to the surjectivity, let $(a,f)\in C_q\,$, so $q(a)=f(0)\,$. There exists $b'\in\Cst0{\mathcal C_F}$ \scth $p_1(b')=a\,$. Then $p(b')=(a,f')\,$. We find $f'(0)=q(a)=f(0)\,$,
  so $f-f'$ belongs to $\ker e_0=I(\Opint001\times F)=\ker p_1\,$. Hence, $f-f'$ can be considered as an element of
  $\Cst0{\mathcal C_F}\,$, and $p_1(f-f')=0\,$. Thus, setting $b=b'+f-f'\,$, we find $p(b)=(a,f')+(0,f-f')=(a,f)\,$. This shows that $p$ is surjective.

  The same argument goes through in the reduced case if $\mathcal G|F$ is topologically amenable, since $\Ropint001\times\mathcal G|F$ is then also amenable.
\end{pf}

The groupoid expression of the mapping cone gives an easy proof of naturality.

\begin{prop}\label{ext-nat}
  Let $\mathcal H$ be another locally compact $\sigma$-compact groupoid, and $G\subset\mathcal H^{(0)}$ a closed invariant subset.
  Assume that $\vphi:\mathcal H\to\mathcal G$ is a strict morphism \scth $\vphi(G)\subset F\,$. Let $\vphi''$ be the restriction of $\vphi$ to $\mathcal H|G\,$, \sth the 
  right square in the following diagram commutes:
  \[
  \xymatrix{%
    (\mathcal G|F)\times\Opint001\ar[r]&\mathcal C_F&\mathcal G\ar[l]&\mathcal G|F\ar[l]\\
    (\mathcal H|G)\times\Opint001\ar[r]\ar[u]^-{\vphi''\times\id_{\Opint001}}&\mathcal C_G\ar@{.>}[u]^-{\psi}&\mathcal H\ar[l]\ar[u]_-{\vphi}
    &\mathcal H|G\ar[l]\ar[u]_-{\vphi''}
  }
  \]
  Here, each of the horizontal arrows is given by an either closed or open inclusion. Then there exists a strict morphism $\psi$ as indicated, which is proper if $\vphi$ is,
  \scth the entire diagram becomes commutative.
\end{prop}

\begin{pf}
  Indeed, simply set $\psi=\vphi\cup_{\mathcal H|G}\bar\vphi$ where $\bar\vphi=\vphi''\times\id_{\Lopint001}\,$. Then the diagram is commutative. The inclusions
  $\xymatrix@1@M+1pt{\mathcal G\ar[r]&\mathcal C_F&\mathcal G|F\ar[l]}$ are closed embeddings: hence, they are proper. If $\vphi$ is proper,
  then so is $\bar\vphi\,$. If $K\subset\mathcal C_F$ is compact, then, identifying subsets of $\mathcal H$ and $\mathcal H|G\times\Lopint001$ with their images
  in $\mathcal C_G\,$, $\psi^{-1}(K)=\vphi^{-1}(K)\cup\bar\vphi^{-1}(K)\,$, which is compact as the union of two compacts. Thus, $\psi$ is proper.
\end{pf}

Let $U=\mathcal G^{(0)}\setminus F$ and $V=\mathcal H^{(0)}\setminus G\,$. By construction, the restriction of $\psi$ to $\mathcal H|V$ is simply $\vphi'\,$,
which sends $\mathcal H|V\to\mathcal G|U\,$. We obtain the following corollary.

\begin{cor}\label{conn-hom-nat-kk}
  Retain the notation of Proposition \ref{ext-nat} and let $\vphi'$ be the restriction of $\vphi$ to $\mathcal H|V\,$.  If $\vphi$ is proper and the groupoids $\mathcal G$ 
  and $\mathcal H$ are amenable and have Haar systems, then the following diagram commutes in $KK\,$:
  \[
  \xymatrix@C+0.5ex{%
    S\Cred0{\mathcal G|F}\ar[r]^-{y\otimes\partial}\ar[d]_-{SKK(\vphi')}&\Cred0{\mathcal G|U}\ar[d]^-{KK(\vphi'')}\\
    S\Cred0{\mathcal H|G}\ar[r]^-{y\otimes\partial}&\Cred0{\mathcal H|V}
  }
  \]
  Here, $S$ denotes suspension and the horizontal maps are the connecting maps in $KK$ theory.
\end{cor}

\begin{pf}
  Applying the cofunctor $KK(-)$ introduced in Section \ref{preliminaries-section}, the following diagram commutes in the category $KK\,$:
  \[
  \xymatrix{%
    S\Cred0{\mathcal G|F}\ar[r]\ar[d]_-{SKK(\vphi')}&C_q=\Cred0{\mathcal C_F}\ar[d]^-{KK(\vphi')}\\
    S\Cred0{\mathcal H|G}\ar[r]&C_p=\Cred0{\mathcal C_G}
  }
  \]
  where $C_q$ and $C_p$ are the mapping cones for 
  \[
  q:\Cred0{\mathcal G}\to\Cred0{\mathcal G|F}\nd p:\Cred0{\mathcal H}\to\Cst0{\mathcal H|G}\ ,
  \]
  and the horizontal maps are natural. As previously observed, the connecting map for $\mathcal G$ is the inverse of the $KK$ equivalence
  $\Cred0{\mathcal G|U}\to C_q\,$, in turn induced by the open inclusion $\mathcal G|U\subset\mathcal G\,$, cf.~\cite[Theorem 19.5.7]{blackadar-ktheory}. Similarly,
  this applies to $\mathcal H\,$. As noted above, $KK(\psi)$ pushes through these equivalences to the arrow $KK(\vphi')\,$. 
\end{pf}

\section{Fibrewise Differentiable Groupoids}

\subsection{Basic Definitions}

  In this section, we extend the concepts of continuous family manifolds and groupoids (of class $\mathcal C^{\infty,0}$), introduced by Paterson 
  \cite{paterson-contfamgrp}, to the case of finite differentiability class $\mathcal C^{q,0}\,$, $q<\infty\,$. This goes through without much ado, and 
  we do not claim any originality in this respect. To our purposes, the interesting point is that the appropriate generalisation of Connes's tangent groupoid can be defined 
  in the finite differentiability class $\mathcal C^{q,0}\,$, $q\sge1\,$. Given Paterson's thorough treatment of the $q=\infty$ case, we need only sketch the elements 
  of the theory for $q<\infty\,$. Continuous family groupoids have been studied in their own right; we refer the interested reader to 
  \cite{paterson-contfamgrp,lauter-monthubert-nistor-pdocfg}.
    
  Let $Y$ be a paracompact topological space, and $A\subset Y\times\reals^n\,$, $B\subset Y\times\reals^m$ be open. Then a continuous and
  fibre-preserving map $f:A\to B$ is said to be of class $\mathcal C^{q,0}\,$, where $q\in\nats\cup\infty\,$, if for any
  $U\times V\subset A$ and $U'\times V'\subset B$ where $U,U'\subset Y$ and $V\subset\reals^n\,$, $V'\subset\reals^m$ are open subsets and
  $f(U\times V)\subset U'\times V'\,$, the map
  \[
  U\to U'\times\Ct[^q]0{V,V'}:y\mapsto f^y=f(y,\var)
  \]
  is well-defined and continuous for the usual Fr\'echet topology on $\Ct[^q]0{U,U'}\,$. For $n=m\,$, the collection of all the locally defined $\mathcal C^{q,0}$ maps with 
  $\mathcal C^{q,0}$ inverse defines a pseudogroup $\Gamma^{q,0}_Y(\reals^n)$ of local homeomorphisms, cf.~\cite{kobayashi-nomizu-vol1}.

  Let $M,Y$ be paracompact locally compact Hausdorff spaces, and let the map $p:M\to Y$ be a continuous open surjection. Then $(M,p)$ is called a manifold of
  class $\mathcal C^{q,0}$ over $Y$ if it has an atlas compatible with $\Gamma_Y^{q,0}(\reals^n)\,$. In this case, each of the
  fibres $M^y=p^{-1}(y)$ is a manifold of class $\mathcal C^q\,$, with dimension constant by definition. It is clear
  how to define the appropriate morphisms to turn the class of $\mathcal C^{q,0}$ manifolds over the fixed base $Y$ into a category.

  Moreover, given a continuous map $f:Z\to Y\,$, any $Y$-manifold $(M,p)$ of class $\mathcal C^{q,0}$ pulls back to a $Z$-manifold $f^*M$ of class
  $\mathcal C^{q,0}\,$. This enables one to define $\mathcal C^{q,0}$ maps between $\mathcal C^{q,0}$ over different bases, and 
  gives rise to a sensible category of $\mathcal C^{q,0}$ manifolds with arbitrary base. 

  Similarly as above, we may define a pseudogroup $\GL_Y^{q,0}(\reals^n,\reals^k)$ of local homeomorphisms of $Y\times\reals^n\times\reals^k$ by 
  considering maps $f$ such that 
  \[
  f(y,a,x)=(y,f_y(a),L_{y,a}(x))\mathtxt{where}L_{y,a}\in\GL(k,\reals)\ .
  \]
  This allows the definition of vector bundles of class $\mathcal C^{q,0}\,$. A trivial but striking consequence is that if $E$ is a topological vector bundle over $Y\,$, 
  and $p:M\to Y$ is a $\mathcal C^{q,0}$ manifold, then $p^*E$ is naturally a $\mathcal C^{q,0}$ vector bundle over $M\,$. 

  Let $q\sge1$ and $(M,p)$ be a manifold over $Y$ of class $\mathcal C^{q,0}\,$. Then we define the fibrewise tangent bundle $TM$ as follows.
  Set-theoretically, $TM$ is the (disjoint) union $TM=\bigcup_{y\in Y}TM^y$ where $M^y=p^{-1}(y)\,$, and the bundle projection
  is $\pi(y,x,\xi)=x\,$. Let $\Parens0{(U_\alpha,\phi_\alpha)}$ be an atlas for $(M,p)\,$, compatible with
  $\Gamma_Y^{p,0}(\reals^n)\,$. Then let $\phi_\alpha^y=\phi_\alpha|(M^y\cap U_\alpha)\,$, and
  \[
  \psi_\alpha:\pi^{-1}(U_\alpha)\to U_\alpha\times\reals^k:(x,\xi)\mapsto\Parens0{\phi_\alpha(x),T_x\phi_\alpha^{p(x)}\xi}\ .
  \]
  (If $y=p(x)$ and $\xi=\dot x(0)$ where $x(t)$ is a $\mathcal C^1$ curve in $M^y\,$, $x(0)=x\,$, then
  $T_x\phi_\alpha^y\xi=\dot z(0)$ where $z(t)=\phi_\alpha\circ x(t)=\phi_\alpha^y\circ x(t)\,$.) Endow $TM$ with the weakest
  topology turning all the $\psi_\alpha$ into homeomorphisms. Then $\Parens0{(\pi^{-1}(U_\alpha),\psi_\alpha)}$ is the
  structure of a vector bundle over $M$ of class $\mathcal C^{q-1,0}$ and rank $k=\dim M^y\,$.

  Assume $(M,p)$ and $(M',p')$ are spaces of class $\mathcal C^{q,0}$ where $q\sge1\,$, and $f:M\to M'$ is a class $\mathcal C^{q,0}$
  morphism. Then we may define a class $\mathcal C^{q-1,0}$ morphism $Tf:TM\to TM'$ in the obvious way, using fibrewise differentiation.

  \begin{defn}\label{def-cq0-grp}
      A groupoid $\mathcal G$ is said to be of class $\mathcal C^{q,0}$ over $\mathcal G^{(0)}$ if
      \begin{enumerate}
      \item $(\mathcal G,r)$ and $(\mathcal G,s)$ are $\mathcal G^{(0)}$-manifolds of class $\mathcal C^{q,0}\,$,
      \item inversion is an isomorphism of class $\mathcal C^{q,0}$ between $(\mathcal G,r)$ and $(\mathcal G,s)\,$, and
      \item considering $\circ:(\mathcal G^{(2)},\pr_1)\to(\mathcal G,r)\,$, $(\circ,r)$ is a morphism of class $\mathcal C^{q,0}\,$.
      \end{enumerate}
      A $\mathcal C^{q,0}$ homomorphism of groupoids is a groupoid homomorphism $f:\mathcal G\to\mathcal H$ between $\mathcal C^{q,0}$
      groupoids $\mathcal G$ and $\mathcal H$ \scth that $(f,f|\mathcal G^{(0)})$ is a $\mathcal C^{q,0}$ morphism for both
      $(\mathcal G,r)\to(\mathcal H,r)$ and $(\mathcal G,s)\to(\mathcal H,s)\,$.
      
      Let $q\sge1$ and $\mathcal G$ be a groupoid of class $\mathcal C^{q,0}\,$. Considering $\mathcal G^{(0)}\subset\mathcal G\,$, we
      may take the restriction $A(\mathcal G)=T\mathcal G|\mathcal G^{(0)}\,$, the so-called Lie algebroid of $\mathcal G\,$, as a vector bundle on $\mathcal G^{(0)}\,$. 
  \end{defn}

\subsection{The Fibrewise Tangent Groupoid}

  Now we are ready to define the fibrewise tangent groupoid of a groupoid $\mathcal G$ of class $\mathcal C^{q,0}\,$.
  Set-theoretically, this is 
  \[
  \mathbb T\mathcal G=A(\mathcal G)\times0\cup\mathcal G\times\Lopint001\ .
  \]
  The unit space is $(\mathbb T\mathcal G)^{(0)}=\mathcal G^{(0)}\times[0,1]\,$, with source and range maps defined by
  \[
  s(x,\xi,0)=(x,0)\ ,\ s(\gamma,\eps)=s(\gamma)\nd r(x,\xi,0)=(x,0)\ ,\ r(\gamma,\eps)=r(\gamma)\ ,
  \]
  and composition given by
  \[
  (x,\xi_1,0)(x,\xi_2,0)=(x,\xi_1+\xi_2,0)\nd (\gamma_1,\eps)(\gamma_2,\eps)=(\gamma_1\gamma_2,\eps)\ .
  \]
  Consider the product topology on $\mathcal G^{(0)}\times[0,1]\,$. The topology of $\mathbb T\mathcal G$ is the weakest for which $r$ and $s$ are
  continuous, as well as the maps $\mathbb Tf:\mathbb T\mathcal G\to\reals\,$, defined by
  \[
  \mathbb Tf(x,\xi,0)=T_xf(\xi)\nd\mathbb Tf(\gamma,\eps)=\frac{f(\gamma)}\eps
  \]
  for any $\mathcal C^{q,0}$ map $f:\mathcal G\to\reals$ which is $0$ on $\mathcal G^{(0)}\,$. 

\begin{prop}
  For any groupoid $\mathcal G$ of class $\mathcal C^{q,0}\,$, $q\sge1\,$, $\mathbb T\mathcal G$ is a groupoid of class $\mathcal C^{q,0}\,$.
\end{prop}

The proof is the same as in the $\mathcal C^{\infty,0}$ case, and we do not repeat it. 

\begin{rem}\ 
\begin{enumerate}\romannum
\item Note that the differentiability class of $\mathbb T\mathcal G$ is the same as for $\mathcal G\,$. This is due to the fact that only 
  differentiability in the fibre direction is considered, and the fibres of $\mathbb T\mathcal G$ over $\mathcal G^{(0)}\times\eps\,$, $\eps>0\,$, are of class $\mathcal C^q\,$, 
  whereas over $\mathcal G^{(0)}\times0\,$, they are linear and hence of class $\mathcal C^\infty\,$. 
\item The prime example of a $\mathcal C^{q,0}$ groupoid for which the tangent groupoid is considered is $\mathcal G=M\times_YM$ where $p:M\to Y$ is a 
  $\mathcal C^{q,0}$ manifold over $Y\,$, $q\sge1\,$. Then $r,s$ are the projections and composition is the same as for the pair groupoid. In this case, 
  $A(\mathcal G)=TM\,$, as is easy to see, so 
  \[
  \mathbb T(M\times_YM)=TM\times0\cup\Parens0{M\times_YM\times\Lopint001}\ ,
  \]
  with the weakest topology that makes source and range continuous, along with the maps $\tilde f:\mathbb T(M\times_YM)\to\reals$ defined for any 
  $f\in\Ct[^{q,0}]0{M,\reals}$ by
  \[
  \tilde f(x,\xi,0)=T_xf(\xi)\nd\tilde f(x_1,x_2,\eps)=\frac{f(x_1)-f(x_2)}\eps\ .
  \]
  Indeed, for any such $f\,$, $h_f:M\times_YM\to\reals\,$, $h_f(x_1,x_2)=f(x_1)-f(x_2)$ is a map of class $\mathcal C^{q,0}\,$. 

  For $Y=\mathrm{pt}\,$, $\mathbb T(M\times_YM)=\mathbb T(M\times M)$ is Connes's tangent groupoid \cite[\S~II.5]{connes-noncommgeom} for the manifold $M\,$. In 
  general, one may view $\mathbb T(M\times_YM)$ as a `family version' of this tangent groupoid. 
\end{enumerate}
\end{rem}

  The following result follows immediately from the triviality of the density bundle $\Abs0\Omega(T^*\mathcal G)\,$, cf.~\cite[proof of Theorem 1]{paterson-contfamgrp}.

\begin{prop}\label{contfamgrp-haarsyst}
  Let $\mathcal G$ be a groupoid of class $\mathcal C^{q,0}\,$, $q\sge1\,$. Then $\mathcal G$ has a Haar system $(\lambda^u)$ which is locally of the form
  $\lambda^u|U=\delta_u\otimes\alpha^u_U\,$, where $\alpha^u_U$ is absolutely continuous to Lebesgue measure on an open subset of $\reals^n\,$.
\end{prop}

\begin{prop}
  Let $\mathcal G$ be a groupoid of class $\mathcal C^{q,0}\,$, $q\sge1\,$. Then $A(\mathcal G)$ is topologically amenable. In particular,
  we have a short exact sequence
  \[
  \xymatrix{%
    0\ar[r]&{\Ct[_0]0{\Lopint001,\Cred0{\mathcal G}}}\ar[r]&{\Cred0{\mathbb T\mathcal G}}\ar[r]^-{e_0}&{\Cred0{A(\mathcal G)}}\ar[r]&0%
  }
  \]
  of reduced groupoid C$^*$-algebras. Here, we denote by 
  \[
  e_t:\Cst0{\mathbb T\mathcal G}\to\Cst0{\mathbb T\mathcal G|\mathcal G^{(0)}\times t}
  \] 
  the maps induced by restriction to the closed invariant subsets 
  \[
  \mathcal G^{(0)}\times t\subset\mathcal G\times[0,1]=(\mathbb T\mathcal G)^{(0)}\ .
  \]
\end{prop}

\begin{pf}
  The groupoid $A(\mathcal G)$ is amenable. Indeed, its isotropy groups $T_x\mathcal G^x$ are Abelian and hence amenable as groups. The principal groupoid associated to
  $A(\mathcal G)$ is the graph of the identity on $\mathcal G^{(0)}\,$, so it is just the space $\mathcal G^{(0)}\,$. The latter is amenable by definition. Then
  \cite[Corollary 5.3.33]{anantharaman-delaroche-renault} gives the measurewise amenability of $A(\mathcal G)\,$; the topological amenability follows from
  \cite[Theorem 3.3.7]{anantharaman-delaroche-renault} (the orbits in $A(\mathcal G)^{(0)}=\mathcal G^{(0)}$ are points). The remaining claims follow from
  section \ref{short-ex-seq}.
\end{pf}

\begin{defn}
  The C$^*$-algebra $\Ct[_0]0{\Lopint001,\Cred0{\mathcal G}}$ being C$^*$-contractible, the map $e_0$ is a $KK$ equivalence, thus inducing an element
  \[
  \tau=e_0^{-1}\otimes e_1^{\vphantom{-1}}\in KK(\Cred0{A(\mathcal G)},\Cred0{\mathcal G})
  \]
  which we call the \emph{Connes--Skandalis map}, cf.~\cite[Definition 3.2]{hilsum-skandalis}. 
\end{defn}

  In fact, such a map can be introduced for any continuous field of C$^*$-algebras over $[0,1]$ which is trivial over
  $\Lopint001\,$. In our $\mathcal C^{1,0}$ groupoid setup, we shall show that the $KK^1$ class $y\otimes\tau$ represents an extension of groupoid C$^*$-algebras.

\section{The Suspended Connes--Skandalis Map}

\subsection{Suspension and Cone on the Tangent Groupoid}\label{suspcone-tangent}

  Let $\mathcal G$ be a groupoid of class $\mathcal C^{1,0}\,$. To prove our index theorem, we shall have to compute
  $y\otimes\tau$ where $\tau$ is the Connes--Skandalis map associated to $\mathbb T\mathcal G\,$, and the element $y\in KK^1(S,\cplxs)$ represents the 
  Wiener--Hopf extension, i.e., represents the standard extension of $\mathcal W_{\reals_{\sge0}}\,$. Whereas $x\otimes\tau$ 
  (where $x=y^{-1}$) is easily evaluated without resorting to groupoid constructions (cf.~\cite[Remark 3.3.2]{hilsum-skandalis}), we shall have to construct an 
  auxiliary groupoid in order to compute the suspended Connes--Skandalis map $y\otimes\tau\,$.

  Recall that $\mathcal W_{\reals_{\sge0}}=(\reals\rtimes\reals)|\reals_{\sge0}\cup(\infty\times\reals)$ is the disjoint union of groupoids, with the topology given as a
  subspace of $[0,\infty]\times\reals\,$. As a topological space, let $\mathbb W\mathcal G=\mathbb T\mathcal G\times_{[0,\infty]}\mathcal W_{\reals_{\sge0}}$ where
  the map $\mathbb T\mathcal G\to[0,\infty]$ is the composition of $r$ (or $s$) with
  \[
  \mathcal G^{(0)}\times[0,1]\to[0,\infty]:
  \begin{cases}
    (x,\eps)\mapsto\frac1\eps-1&\eps>0\ ,\\
    (x,0)\mapsto\infty &\text{otherwise},
  \end{cases}
  \]
  and $\mathcal W_{\reals_{\sge0}}\to[0,\infty]$ is the range projection. Define groupoid operations on $\mathbb W\mathcal G$ as follows:
  \begin{gather*}
    r(\gamma,r_1,r_2-r_1)=(r(\gamma),r_1)\ ,\ s(\gamma,r_1,r_2-r_1)=(s(\gamma),r_2)\ ,\\
    r(x,\xi,\infty,r)=(x,\infty)=s(x,\xi,\infty,r)\ ,    
  \end{gather*}
  and
  \begin{gather*}
    (\gamma_1,r_1,r_2-r_1)(\gamma_2,r_2,r_3-r_2)=(\gamma_1\gamma_2,r_1,r_3-r_1)\ ,\\
    (x,\xi_1,\infty,r_1)(x,\xi_2,\infty,r_2)=(x,\xi_1+\xi_2,\infty,r_1+r_2)\ .
  \end{gather*}

\begin{prop}\label{wg-haarsyst}
  Given a class $\mathcal C^{1,0}$ groupoid $\mathcal G\,$, the space $\mathbb W\mathcal G$ is a locally compact groupoid \scth
  $\mathbb W\mathcal G^{(0)}=\mathcal G^{(0)}\times[0,\infty]\,$. The subset $F=\mathcal G^{(0)}\times\infty$ is closed and invariant, and we have
  \[
  \mathbb W\mathcal G|F=A(\mathcal G)\times\reals\nd\mathbb W\mathcal G|U=\mathcal G\times(\reals\rtimes\reals)|\reals_{\sge0}\mathtxt{for}
  U=\mathbb W\mathcal G^{(0)}\setminus F\ .
  \]
  Moreover, $\mathbb W\mathcal G$ carries a natural Haar system which may be chosen to induce on $\mathcal G$ any given Haar system associated to a positive section of
  the density bundle $\Abs0\Omega(T^*\mathcal G)\,$.
\end{prop}

\begin{pf}
  It is clear that $\mathbb W\mathcal G$ is a locally compact space, and it is also evidently a groupoid. We have the following commutative diagram:
  \[
  \xymatrix{%
    {\mathbb T\mathcal G^{(2)}\times\mathcal W_{\reals_{\sge0}}^{(2)}}\ar[r]^-{\circ}&{\mathbb T\mathcal G\times\mathcal W_{\reals_{\sge0}}}\\
    {\mathbb W\mathcal G^{(2)}}\ar[r]_-{\circ}\ar[u]&\mathbb W\mathcal G\ar[u]}
  \]
  where the vertical maps are self-evident, and the rightmost of these is a closed embedding. Thus, composition is continuous, and along the same lines, the continuity of the
  inverse is established. The projections $r$ and $s$ are continuous.

  As to the existence of Haar systems, $\mathcal W=\,\Lopint0{-\infty}\infty\rtimes\reals$ is a $\mathcal C^{1,0}$ groupoid; in fact, the direct product
  $\mathcal W=\,\Lopint0{-\infty}\infty\times\reals$ of spaces is a $\mathcal C^{1,0}$ manifold over $\Lopint0{-\infty}\infty\,$, and the operations are fibrewise
  those of the Lie group $\reals\,$, independent of the fibre. Moreover, $\mathcal W_{\reals_{\sge0}}$ is the restriction of $\mathcal W$ to the non-invariant subset
  $[0,\infty]$ of $\Lopint0{-\infty}\infty\,$.

  Similarly as for $\mathbb T\mathcal G\,$, we may define a $\mathcal C^{1,0}$ groupoid $\mathcal T=A(\mathcal G)\times0\cup\mathcal G\times\,\Opint00\infty\,$, 
  by replacing $[0,1]$ in the definition of the tangent groupoid by $\Ropint00\infty\,$. Choosing a homeomorphism
  $\phi:\Ropint00\infty\,\to\,\Lopint0{-\infty}\infty$ which coincides on $[0,1]$ with $\eps\mapsto\tfrac1\eps-1\,$, we obtain a $\mathcal C^{1,0}$ groupoid
  $\mathcal H=\mathcal T\times_{\Lopint0{-\infty}\infty}\mathcal W$ \scth $\mathbb W\mathcal G=\mathcal H|(\mathcal G^{(0)}\times[0,\infty])\,$. In fact, if
  $f:\mathcal G\to\reals^n$ is a $\mathcal C^{1,0}$ map \scth $(r,f)$ is a local chart, we may define $\psi_f:\mathcal H\to\reals^{n+1}$ by
  \[
  \psi_f(\tau)=
  \begin{cases}
    \Parens1{\tfrac{f(\gamma)}{\phi^{-1}(r_1)},r_2-r_1}&\tau=(\gamma,r_1,r_2-r_1)\,,\,r_1<\infty\ ,\\
    \Parens0{T_xf(\xi),r}&\tau=(x,\xi,\infty,r)\ .
  \end{cases}
  \]
  Then $(r,\psi_f)$ is a local chart for $\mathcal H\,$. Now, $\mathcal H$ has a Haar system given by a positive section of the density bundle, unique up to
  multiplication by such a density. Thus we may assume that this Haar system induces on $\mathcal G$ the given Haar system associated to the choice of a positive section
  of $\Abs0\Omega(T^*\mathcal G)\,$.

  If $\lambda^{x,t}\,$, $(x,t)\in\mathcal G^{(0)}\times\,\Lopint0{-\infty}\infty\,$, is a Haar system of $\mathcal H\,$, define an invariant
  system of positive Radon measures by $\mu^{x,t}=\lambda^{x,t}|\mathbb W\mathcal G^{x,t}\,$. Since $\mathbb W\mathcal G^{x,t}$ has dense
  interior in $\mathcal H^{x,t}\,$, the measures $\mu^{x,t}$ satisfy the support condition. The maps $x\mapsto\lambda^{x,t}\,$, for $t\in[0,\infty]\,$, are
  equicontinuous. Hence, the same is true for $x\mapsto\mu^{x,t}\,$. Since for fixed $x\,$, the characteristic functions of the interiors of
  $\mathbb W\mathcal G^{x,t}$ depend continuously in the topology of simple convergence on $t\,$, we find that $\mu^{x,t}$ satisfies the continuity axiom.
  The statement about the invariant subsets and the corresponding restricted groupoids is immediate.
\end{pf}

\begin{cor}\label{wg-ext}
  There is a short exact sequence
  \[
  \xymatrix@C-16pt{%
    0\ar[r]&\Cred0{\mathcal G}\otimes\mathbb K\ar[r]&\Cred0{\mathbb W\mathcal G}\ar[r]%
    &S\Cred0{A(\mathcal G)}=\Cred0{A(\mathcal G)}\otimes\Cred0\reals\ar[r]&0}\ .\tag{$**$}
  \]
\end{cor}

\begin{pf}
  We need only remark that $A(\mathcal G)\times\reals$ is an amenable groupoid, and that $(\reals\rtimes\reals)|\reals_{\sge0}$ is isomorphic to the pair groupoid
  $\reals_{\sge0}\times\reals_{\sge0}\,$, whose reduced C$^*$-algebra is $\mathbb K\,$.
\end{pf}

  To see that $\tau$ `interpolates' between the Wiener--Hopf extension and the one constructed above, we need to construct the `cone' $\mathbb C\mathcal G$ over the
  tangent groupoid. This is the content of the following proposition.

\begin{prop}\label{cone-constr}
  Let $\mathcal G$ be a groupoid of class $\mathcal C^{1,0}\,$. There exists a locally compact groupoid $\mathbb C\mathcal G$
  over the `triangle'
  \[
  \mathbb C\mathcal G^{(0)}=\mathcal G^{(0)}\times\Delta\mathtxt{where}\Delta=[0,1]\times[0,\infty]/[0,1]\times\infty\ ,
  \]
  \scth $U=[0,1]\times\mathcal G^{(0)}\times\Lopint00\infty$ is an open invariant subset of $\mathbb C\mathcal G^{(0)}\,$, 
  \[
  \mathbb C\mathcal G|F=A(\mathcal G)\times\reals\nd\mathbb C\mathcal G|U=\mathbb T\mathcal G\times(\reals\rtimes\reals)|\reals_{\sge0}\mathtxt{where}
  F=\mathbb C\mathcal G^{(0)}\setminus U\ .
  \]
  In addition, $\mathbb C\mathcal G$ carries a Haar system which induces on $\mathbb T\mathcal G$ a Haar system given by the choice of a positive section of the latter
  groupoid's density bundle.
\end{prop}

\begin{pf}
  Let $\mathcal H$ be the $\mathcal C^{1,0}$ groupoid over $\mathcal H^{(0)}=\mathcal G^{(0)}\times\Lopint0{-\infty}\infty$ constructed in the proof of
  Proposition \ref{wg-haarsyst}, \scth $\mathbb W\mathcal G=\mathcal H|\Parens1{\mathcal G^{(0)}\times[0,\infty]}\,$. We construct the `partial' tangent groupoid
  $\mathcal T\mathcal H=0\times A(\mathcal G)\times\mathcal W\,\cup\,\Lopint001\times\mathcal H\,$, 
  the disjoint union of groupoids. We endow this set with the initial topology with respect to $r\,$, $s\,$, and the maps
  $\vrho_f:\mathcal{TH}\to\reals^{n+1}$ defined for $\mathcal C^{1,0}$ charts $(r,f)$ of $\mathcal G\,$, $f:\mathcal G\to\reals^n\,$, as follows:
  \[
  \vrho_f(\tau)=
  \begin{cases}
    \Parens1{\tfrac{f(\gamma)}{\eps+\phi^{-1}(r_1)},r_2-r_1}&\tau=(\eps,\gamma,r_1,r_2-r_1)\,,\,\eps>0\,,\,r_1<\infty\ ,\\
    \Parens0{T_xf(\xi),r_2-r_1}&\tau=(0,x,\xi,r_1,r_2-r_1)\,,\,r_1<\infty\ ,\\
    \Parens0{T_xf(\xi),r}&\tau=(\eps,x,\xi,\infty,r)\,,\,\eps\in[0,1]\ .
  \end{cases}
  \]
  Then, for any such $f\,$, $(r,\vrho_f)$ is a local $\mathcal C^{1,0}$ chart for $\mathcal{TH}\,$, turning the latter into a $\mathcal C^{1,0}$ groupoid.

  Define on the unit space of $\mathcal{TH}\,$, $\mathcal{TH}^{(0)}=[0,1]\times\mathcal G^{(0)}\times\Lopint0{-\infty}\infty\,$, the following equivalence relation: 
  \begin{align*}
    (\eps_1,x_1,r_1)&\sim(\eps_2,x_2,r_2)\\ &\equiva\ x_1=x_2\nda\Parens1{\min(r_1,r_2)<\infty\ \imp (\eps_1,r_1)=(\eps_2,r_2)}\ ,    
  \end{align*}
  and denote its graph by $S\,$. Then $S$ is a subgroupoid of the pair groupoid on $\mathcal{TH}^{(0)}\,$, and it acts on $\mathcal{TH}$ by
  \[
  (s,t).\gamma=
  \begin{cases}
    \gamma &t=r(\gamma)\not\in[0,1]\times\mathcal G^{(0)}\times\infty\ ;\\
    (\eps_1,x,\xi,\infty,r)&
    \left\{\begin{aligned}[c]
      &s=(\eps_1,x,\infty)\,,\,t=(\eps_2,x,\infty)\,,\\
      &\gamma=(\eps_2,x,\xi ,\infty,r)\ .
    \end{aligned}\right.
  \end{cases}
  \]
  Thus, $S$ fixes $\gamma$ whenever $r(\gamma)\not\in[0,1]\times\mathcal G^{(0)}\times\infty\,$, and on
  \[
  r^{-1}([0,1]\times\mathcal G^{(0)}\times\infty)=[0,1]\times A(\mathcal G)\times\infty\times\reals\ ,
  \]
  the action of $S$ fixes all but the first component, and there, it is the usual action of the pair groupoid over $[0,1]\,$.

  Let $R$ denote the graph of the equivalence relation on $\mathcal{TH}$ defined by the action of $S\,$. Then $R$ is a closed subset of $\mathcal{TH}\times\mathcal{TH}\,$,
  and the equivalence classes of $R$ are compact. Therefore, $\mathcal{CH}=\mathcal{TH}/R$ is a locally compact space, and the associated canonical projection
  $\pi:\mathcal{TH}\to\mathcal{CH}$ is proper, by \cite[Chapter I, \S~10.4, Proposition~9]{bourbaki-topology1}. Moreover, the charts $\vrho_f$ are invariant for the action of $S\,$,
  and hence drop to $\mathcal{CH}\,$, thereby turning this space into a $\mathcal C^{1,0}$ manifold over the `triangle' $\mathcal{CH}^{(0)}=\mathcal{TH}^{(0)}/S\,$. In fact,
  the operations of $\mathcal{CH}$ commute with the action of $S\,$, and since they are compatible with the charts $\vrho_f\,$, $\mathcal{CH}$ turns into a
  $\mathcal C^{1,0}$ groupoid.

  In particular, $\mathcal{CH}$ has a Haar system induced by the choice of a positive density. It restricts to a Haar system for
  \[
  \mathbb C\mathcal G=\mathcal{CH}|(\mathcal G^{(0)}\times\Delta)\mathtxt{where}\Delta=[0,1]\times[0,\infty]/[0,1]\times\infty\ ,
  \]
  by the same argument as in the proof of Proposition \ref{wg-haarsyst}.

  Let $U'=[0,1]\times\mathcal G^{(0)}\times\reals\,$. Because the quotient map $[0,1]\times[0,\infty]\to\Delta$ restricts to a homeomorphism on 
  $[0,1]\times\Ropint00\infty\,$, the restriction of $\pi$ to $\mathcal{TH}|U'$ has local sections and is injective. Hence, it is a homeomorphism onto
  its image. Moreover, $\mathcal{TH}|U'$ is $R$-saturated, so the image is open in $\mathcal{CH}\,$. Obviously, $\mathbb C\mathcal G\cap\pi(U')=U\,$, where 
  $U$ is the image of $[0,1]\times\mathcal G^{(0)}\times\Ropint00\infty$ in $\mathbb C\mathcal G^{(0)}\,$. Thus, the intersection of $\pi(\mathcal{TH}|U')$ with 
  $\mathbb C\mathcal G$ is equal to $\mathbb C\mathcal G|U\,$, and 
  \[
  \mathbb C\mathcal G|U\cong\mathcal{TH}|U=\mathbb T\mathcal G\times(\reals\rtimes\reals)|\reals_{\sge0}\ ,
  \]
  where this isomorphism is the restriction of an isomorphism of $\mathcal C^{1,0}$ groupoids.

  Now, $F=\mathbb C\mathcal G^{(0)}\setminus U$ equals the complement of $U'$ in $\mathcal{TH}^{(0)}\,$. Since the action of $S$ on $\mathcal{TH}|F$ identifies
  with the standard action of the pair groupoid on $[0,1]\,$, we have
  \[
  \mathbb C\mathcal G|F=\mathcal{TH}|F\cong[0,1]\times A(\mathcal G)\times\reals/[0,1]\times[0,1]\cong A(\mathcal G)\times\reals\ ,
  \]
  where the latter is an isomorphism of $\mathcal C^{1,0}$ groupoids which is fibrewise the identity. This proves our assertion.
\end{pf}

\begin{cor}\label{cone-extension}
  There is a short exact sequence
  \[
  \xymatrix@C-16pt{%
    0\ar[r]&\Cred0{\mathbb T\mathcal G}\otimes\mathbb K\ar[r]&\Cred0{\mathbb C\mathcal G}\ar[r]%
    &S\Cred0{A(\mathcal G)}=\Cred0{A(\mathcal G)}\otimes\Cred0\reals\ar[r]&0}\ .
  \]
\end{cor}

\begin{lem}\label{cone-amenable}
  If $\mathcal G$ is topologically amenable, then so are $\mathbb W\mathcal G$ and $\mathbb C\mathcal G\,$.
\end{lem}

\begin{pf}
  Retain the previous notation. The obvious continuous surjection $p:\mathcal{TH}\to[0,1]$ factors through $r$ and $s\,$. 
  Then $p$ is open when restricted to $\mathcal{TH}^{(0)}=[0,1]\times\Lopint0{-\infty}\infty\times\mathcal G^{(0)}\,$, and $r$ and $s$ are open since 
  $\mathcal{TH}$ carries a Haar system. Thus, $p$ is open, and defines continuous field of groupoids in the sense of \cite[Definition 5.2]{landsman-ramazan}. Hence, 
  $\mathcal{TH}$ is topologically amenable if this is the case for
  the fibres of $p\,$, by \cite[Corollary 5.6]{landsman-ramazan}. The fibre at $0$ is $A(\mathcal G)\times\mathcal W_{\reals_{\sge0}}\,$, which is always amenable. The
  fibre at $\eps>0$ is isomorphic to $\mathcal H\,$, so $\mathcal{TH}$ is amenable if $\mathcal H$ is. By the same argument, $\mathcal T$ is amenable if
  $\mathcal G$ is. But $\mathcal H$ is the fibred product of $\mathcal T$ and $\mathcal W_{\reals_{\sge0}}\,$, so it is amenable if $\mathcal G$ is. So, in this case,
  both $\mathcal{TH}$ and $\mathbb W\mathcal G\,$, as a restriction of $\mathcal H\,$, are amenable. Since $\pi:\mathcal{TH}\to\mathcal{CH}$ is proper, the
  amenability of $\mathcal{CH}\,$, and hence, of its restriction $\mathbb C\mathcal G\,$, also follow.
\end{pf}

\subsection{Computation of the Suspended Connes--Skandalis Map}\label{susp-conneskand}
  Now, we can finally compute $y\otimes\tau\,$, as announced.

\begin{thm}\label{ytensortau}
  Let $\mathcal G$ be a topologically amenable, locally compact groupoid of class $\mathcal C^{1,0}\,$. Let 
  $\tau\in KK(\Cred0{A(\mathcal G)},\Cred0{\mathcal G})$ denote the Connes--Skandalis map for the tangent groupoid $\mathbb T\mathcal G\,$, and $y\in KK^1(S,\cplxs)$
  represent the Wiener--Hopf extension. Then we have
  \[
  y\otimes\tau=\partial\ ,
  \]
  where $\partial\in KK^1(\Cred0{A(\mathcal G)}\otimes\Cred0\reals,\Cred0{\mathcal G})$ represents the extension $(**)$.
\end{thm}

\begin{pf}
  We retain the notations from the proof of Proposition \ref{cone-constr}. We have the commutative diagram of strict homomorphisms
  \[
  \xymatrix@M+1pt{%
    \mathbb W\mathcal G\ar[d]\ar[rr]^-{\phi_1}&&\mathbb C\mathcal G\ar[d]\\
    \mathcal H\ar[r]_-{(1,\id)}&\mathcal T\mathcal H\ar[r]_-{\pi}&\mathcal{CH}}
  \]
  The vertical arrows are closed embeddings. The quotient map $\pi$ is proper by the proof of Proposition \ref{cone-constr}. The restriction of $\pi$ to $1\times\mathcal H$ is
  injective, and thus, a closed embedding. Hence, the strict homomorphism $\phi_1$ induced in the above diagram is a closed embedding, in particular, proper.

  It is easy to compute that the following diagram of strict homomorphisms is commutative:
  \[
  \xymatrix@M+1pt@C+8ex{%
    \mathcal G\times(\reals\rtimes\reals)|\reals_{\sge0}\ar[d]_-{i_1\times\id}\ar[r]&\mathbb W\mathcal G\ar[d]^-{\phi_1}
    &A(\mathcal G)\times\reals\ar[l]\ar@{=}[d]\\
    \mathbb T\mathcal G\times(\reals\rtimes\reals)|\reals_{\sge0}\ar[r]&\mathbb C\mathcal G&A(\mathcal G)\times\reals\ar[l]}
  \]
  Here $i_1:\mathcal G\to\mathbb T\mathcal G$ is the inclusion at the fibre over $1\,$, and thus induces the $*$-homomorphism
  $e_1:\Cred0{\mathbb T\mathcal G}\to\Cred0{\mathcal G}\,$. The groupoids involved being amenable by Lemma \ref{cone-amenable}, we may apply
  Corollary \ref{conn-hom-nat-kk} to obtain $\partial_{\mathbb C\mathcal G}\otimes e_1=\partial\,$, where $\partial_{\mathbb C\mathcal G}$ represents the extension
  from Corollary \ref{cone-extension}.

  Similarly, we have a commutative diagram of strict homomorphisms:
  \[
  \xymatrix@M+1pt{%
    A(\mathcal G)\times\mathcal W_{\reals_{\sge0}}\ar[rr]^-{\phi_0}\ar[d]&&\mathbb C\mathcal G\ar[d]\\
    A(\mathcal G)\times\mathcal W\ar[r]_-{(0,\id)}&\mathcal{TH}\ar[r]_-{\pi}&\mathcal{CH}}
  \]
  Again, the vertical arrows are closed embeddings, as is the restriction of $\pi$ to $A(\mathcal G)\times\mathcal W\,$. Thus, $\phi_0\,$, induced by the above diagram,
  is a proper strict homomorphism. We have a commutative diagram
  \[
  \xymatrix@M+1pt@C+8ex{%
    A(\mathcal G)\times(\reals\rtimes\reals)|\reals_{\sge0}\ar[d]_-{i_0\times\id}\ar[r]&A(\mathcal G)\times\mathcal W_{\reals_{\sge0}}\ar[d]^-{\phi_0}
    &A(\mathcal G)\times\reals\ar[l]\ar@{=}[d]\\
    \mathbb T\mathcal G\times(\reals\rtimes\reals)|\reals_{\sge0}\ar[r]&\mathbb C\mathcal G&A(\mathcal G)\times\reals\ar[l]}
  \]
  Here, $i_0:A(\mathcal G)\to\mathbb T\mathcal G$ is the inclusion at the fibre over $0\,$, and induces the $*$-homomorphism and $KK$ equivalence
  $e_0:\Cred0{A(\mathcal G)}\to\Cred0{\mathbb T\mathcal G}\,$. The upper line induces an extension which is represented by
  \[
  \id\otimes\,y\in KK^1(\Cred0{A(\mathcal G)}\otimes\Cred0\reals,\Cred0{A(\mathcal G)})\ .
  \]
  Applying Corollary \ref{conn-hom-nat-kk} entails $\partial_{\mathbb C\mathcal G}\otimes e_0=\id\otimes\,y\,$. Hence,
  \[
  y\otimes\tau=(\id\otimes\,y)\otimes\tau=\partial_{\mathbb C\mathcal G}\otimes e_1=\partial\ ,
  \]
  which was our claim.
\end{pf}

\section{Topological Expression of the Connes--Skandalis Map}\label{topexpr-connesskand}

\subsection{Naturality of Classifying Spaces}\label{c10-grp-class}
  In order to compute our index in topological terms, we shall be particularly interested in the Connes--Skandalis map in the following situation: Fix a 
  manifold $p:M\to Y$ of class $\mathcal C^{1,0}$ over the locally compact, paracompact space $Y\,$, and assume that 
  \emph{$p$ is closed}. Consider the category whose objects are groupoids $\mathcal G$ of class $\mathcal C^{1,0}$ over $\mathcal G^{(0)}=M$ 
  (cf.~Definition \ref{def-cq0-grp}) and whose arrows are the (strict) groupoid morphisms of class $\mathcal C^{1,0}\,$. In this category, the product of $\mathcal G$ with 
  $M\times_YM$ is 
  \[
  \mathcal G\times_M(M\times_YM)=\Set1{(\gamma,x_1,x_2)}{s(\gamma)=x_1\,,\,p(x_1)=p(x_2)}\ ,
  \]
  with structure maps and groupoid composition given by 
  \[
  s(\gamma,x_1,x_2)=x_1\ ,\ r(\gamma,x_1,x_2)=x_2\ ,\ (\gamma_1,x_1,x_2)(\gamma_2,x_2,x_3)=(\gamma_1\gamma_2,x_1,x_3)\ .
  \]
  The Lie algebroid of $\mathcal G\times_M(M\times_YM)$ is $A(\mathcal G)\oplus TM\,$, the direct sum of vector bundles over $M\,$, where $TM$ is the
  fibrewise tangent bundle of $M\,$.

  In particular, for any topological vector bundle $E\to Y\,$, the pullback $p^*E$ is a class $\mathcal C^{1,0}$ groupoid $\mathcal G$ over $M$
  whose Lie algebroid is $p^*E\oplus TM\,$. It is instructive to note that the strict homomorphism
  \[
  p^*E\times_M(M\times_YM)\to E:(x_1,\xi,x_2)\mapsto(p(x_1),\xi)
  \]
  is a Morita equivalence, although we shall not use this fact directly. From the previous sections, we have a Connes--Skandalis map
  \[
  \tau\in KK(\Cred0{p^*E\oplus TM},\Cred0{p^*E\times_M(M\times_YM)})\ .
  \]
  In fact, this is a version `with coefficients' of the topological family index in sense of Atiyah--Singer \cite[Section 3]{as-ind4}. This fact is well-known in the 
  coefficient-free case (i.e.~$E=0$), although it is does not appear to be as well recorded in the literature, at least in this generality. The parameter- and coefficient-free 
  version can be found in Connes's book \cite{connes-noncommgeom} or in \cite{landsman-quant}. We give a brief discussion of this fact, which also explains 
  what we mean by a `topological family index'.

  The basic idea is Connes's realisation (explained in \cite{landsman-quant}) that the groupoids considered above act properly and freely on $\reals^n$ (where the action 
  is obtained along the lines of Atiyah--Singer's construction), so that the Fourier transform can be used to write down $K$-theory 
  isomorphisms between these groupoids and certain topological spaces explicitly. The details will become clear presently.

  Let $\mathcal G$ be a locally compact groupoid with open range map, and let $h:\mathcal G\to V$ be a strict homomorphism, where $V$ is some
  finite-dimensional real inner product space. Define the space $E(\mathcal G,h)=\mathcal G^{(0)}\times V\,$. Then $\mathcal G$ acts on $E(\mathcal G,h)$ from the right by
  \[
  (r(\gamma),v)\gamma=(s(\gamma),v+h(\gamma))\mathfa(r(\gamma),v,\gamma)\in E(\mathcal G,h)\times_{\mathcal G^{(0)}}\mathcal G\ .
  \]
  This action gives rise to the map
  \[
  E(\mathcal G,h)\times_{\mathcal G^{(0)}}\mathcal G\to E(\mathcal G,h)\times E(\mathcal G,h):(r(\gamma),v,\gamma)\mapsto(s(\gamma),v+h(\gamma),r(\gamma),v)
  \]
  The action is free (resp.~proper) if and only if $(r,h,s):\mathcal G\to\mathcal G^{(0)}\times V\times\mathcal G^{(0)}$ is injective (resp.~proper).

  Now, consider $B(\mathcal G,h)=E(\mathcal G,h)/\mathcal G\,$. This is a locally compact space, and as such, a locally compact (cotrivial) groupoid. Its (trivial)
  action on $E(\mathcal G,h)$ is proper and free, and the quotient by this action is $E(\mathcal G,h)$. Since
  the range map of $\mathcal G$ is open, the canonical projection $\pi_h:E(\mathcal G,h)\to B(\mathcal G,h)$ is open and surjective, by \cite[Lemma 2.30]{tu-nonhausgrp}.
  The following easy lemma characterises when $E(\mathcal G,h)$ is a Morita equivalence of $B(\mathcal G,h)$ and $E(\mathcal G,h)\rtimes\mathcal G\,$.

\begin{lem}
Let $\mathcal G$ be a locally compact groupoid with an open range map. Then the following statements are equivalent.
\begin{enumerate}\romannum
\item The groupoid $\mathcal G$ is principal.
\item The space $\mathcal G^{(0)}$ defines a Morita equivalence $\mathcal G\to\mathcal G^{(0)}/\mathcal G\,$.
\item The canonical projection $\pi:\mathcal G\to\mathcal G^{(0)}/\mathcal G$ is proper as a generalised morphism.
\item The space $\mathcal G^{(0)}$ defines a generalised morphism.
\end{enumerate}
In this case, $\pi$ is the inverse of $\mathcal G^{(0)}$ as a generalised morphism.
\end{lem}


  A simple consequence of the lemma is the following naturality of $E(\mathcal G,h)\,$: Let $(\mathcal G,h)$ and $(\mathcal G',h')$ be given, where the
  groupoids $\mathcal G$ and $\mathcal G'$ are locally compact with open range maps, and $h\,$, $h'$ are strict homomorphisms \scth $(r,h,s)$ and $(r,h',s)$ are
  injective and proper. Then a \emph{morphism of pairs} $\vphi:(\mathcal G,h)\to(\mathcal G',h')$ is a strict homomorphism
  $\vphi:\mathcal G\to\mathcal G'$ \scth $h'\circ\vphi=h\,$. Such a morphism of pairs $\vphi$ gives rise to the continuous map
  \[
  E(\vphi):E(\mathcal G,h)\to E(\mathcal G',h'):(x,v)\mapsto(\vphi(x),v)
  \]
  which intertwines the actions of $\mathcal G$ and $\mathcal G'\,$:
  \[
  E(\vphi)\Parens0{(r(\gamma),v)\gamma}=\Parens0{r(\vphi(\gamma)),v}\vphi(\gamma)\mathfa\gamma\in\mathcal G\,,\,v\in V\ .
  \]
  Hence, $E(\vphi)$ induces a map $B(\vphi):B(\mathcal G,h)\to B(\mathcal G',h')\,$. Moreover,
  \[
  E(\vphi)\times\vphi:E(\mathcal G,h)\rtimes\mathcal G\to E(\mathcal G',h')\rtimes\mathcal G'
  \]
  is a strict homomorphism inducing $B(\vphi)$ under the straightforward identification $B(\mathcal G,h)=E(\mathcal G,h)/E(\mathcal G,h)\rtimes\mathcal G\,$. 
  If $\pi_h$ and $\pi_{h'}$ are the canonical projections, we thus have $\pi_{h'}\circ(E(\vphi)\times\vphi)=B(\vphi)\circ\pi_h\,$. By the lemma, we may take inverses, so the
  following diagram of generalised morphisms commutes:
  \[
  \xymatrix@C+20pt{%
    E(\mathcal G,h)\rtimes\mathcal G\ar[r]^-{E(\vphi)\times\vphi}&E(\mathcal G',h')\rtimes\mathcal G'\\
    B(\mathcal G,h)\ar[u]^-{E(\mathcal G,h)}\ar[r]_-{B(\vphi)}&B(\mathcal G',h')\ar[u]_-{E(\mathcal G',h')}}
  \]
  Since the vertical arrows are Morita equivalences, the horizontal arrows are always simultaneously proper as generalised morphisms. For instance, this is the case if $\vphi$ is
  proper as a strict morphism, since the same is then true of $E(\vphi)\times\vphi\,$. However, in this case, $B(\vphi)$ need not be proper as a strict morphism. For the sake of 
  brevity, we shall write $B(\vphi)^*$ instead of $KK(B(\vphi))$ for the induced $KK$ cycle, even if $B(\vphi)$ is only proper as a generalised morphism.

  For any strict homomorphism $h:\mathcal G\to V$ from the locally compact groupoid with Haar system $\mathcal G$ to the finite-dimensional real inner product space 
  $V\,$, an action of $V$ on $\Cred0{\mathcal G}$ is given by
  \[
  \alpha_h(v)(\vphi)(\gamma)=e^{2\pi i\Rscp0v{h(\gamma)}}\cdot\vphi(\gamma)\mathfa v\in V\,,\,\vphi\in\Cc0{\mathcal G}\,,\,\gamma\in\mathcal G\ .
  \]
  Moreover, we have a $*$-isomorphism $\mathcal F_h:\Cred0{\mathcal G}\rtimes_{\alpha_h}V\to\Cred0{E(\mathcal G,h)\rtimes\mathcal G}\,$, 
  \[
  \mathcal F_h(\vphi)(r(\gamma),v,\gamma)=\int_Ve^{-2\pi i\Rscp0vw}\vphi(\gamma,w)\,dw\mathfa\vphi\in\Cc0{\mathcal G\times V}\ ,
  \]
  as is well-known (cf.~\cite[Proposition 7]{connes-noncommgeom}).

\begin{lem}
  The isomorphism $\mathcal F_h$ is natural in the following sense: Given a morphism of pairs $\vphi:(\mathcal G,h)\to(\mathcal G',h')$ \scth $\vphi$ is proper and
  preserves Haar systems, the following diagram of $*$-morphisms commutes:
  \[
  \xymatrix@C+18pt{%
    \Cred0{\mathcal G'}\rtimes_{\alpha_{h'}}V\ar[r]^-{\vphi^*\otimes\,\id}\ar[d]_-{\mathcal F_{h'}}&
    \Cred0{\mathcal G}\rtimes_{\alpha_h}V\ar[d]^-{\mathcal F_h}\\
    \Cred0{E(\mathcal G',h')\rtimes\mathcal G'}\ar[r]_-{(E(\vphi)\times\vphi)^*}&
    \Cred0{E(\mathcal G,h)\rtimes\mathcal G}}
  \]
\end{lem}


  Now, fix a real inner product space $V$ with $\dim_\reals V=2n\,$. Consider the category $\mathrm{Grp}/V$ whose objects are pairs $(\mathcal G,h)$ where $\mathcal G$
  is a locally compact groupoid with a Haar system and $h:\mathcal G\to V$ is a strict homomorphism \scth $(r,h,s)$ is injective and proper, and whose arrows are
  proper strict homomorphisms $\vphi:\mathcal G\to\mathcal G'$ preserving Haar systems and satisfying $h=h'\circ\vphi\,$. Then we have two cofunctors,
  $KK$ and $KK\circ B\,$, $\mathrm{Grp}/V\to KK\,$, given by
  \[
  KK(\mathcal G,h)=\Cred0{\mathcal G}\ ,\ KK\Parens1{\vphi:(\mathcal G,h)\to(\mathcal G',h')}=\vphi^*=KK(\vphi)\ ,
  \]
  and
  \begin{gather*}
  (KK\circ B)(\mathcal G,h)=\Ct[_0]0{B(\mathcal G,h)}\ ,\\
  (KK\circ B)\Parens1{\vphi:(\mathcal G,h)\to(\mathcal G',h')}=B(\vphi)^*=KK(B(\vphi))\ .    
  \end{gather*}
  By naturality of Thom isomorphisms \cite[Proposition 19.3.5]{blackadar-ktheory}, we have a natural isomorphism of functors $\sigma:KK\imp KK\circ B$ given by
  \[
  \sigma_{\mathcal G,h}=t_h\otimes\mathcal F_h^*KK(E(\mathcal G,h))\mathtxt{where}t_h=t_{\alpha_h}\in
  KK(\Cred0{\mathcal G},\Cred0{\mathcal G}\rtimes_{\alpha_h} V)
  \]
  is the Thom element for the action $\alpha_h\,$.

\subsection{The Connes--Skandalis Map is a Topological Family Index}

We make the following basic observation, which is essentially a minor reformulation of \cite[p.~497]{as-ind1}. Let $Z\to Y$ be a class $\mathcal C^{1,0}$ family, and 
$\pi_V:V\to Z$ a $\mathcal C^{1,0}$ vector bundle. Choose a $\mathcal C^{1,0}$ embedding $i_V:V\to Y\times\reals^k$ compatible with a  
$\mathcal C^{1,0}$ embedding of $Z\,$. Form the normal bundle $NV$ of $V\,$, which is a $\mathcal C^{1,0}$ 
vector bundle over $Z\,$. Then the trivial bundle $NV\times\reals^k$ is homeomorphic to $\pi_V^*(NV\otimes\cplxs)\,$, and hence, a complex vector bundle over $V\,$. In 
particular, there is a Thom isomorphism $\vphi_V:K^*_c(V)\to K_c^*(NV\times\reals^k)\,$. 

Now, let $p:M\to Y$ be a $\mathcal C^{1,0}$ family, $\pi_E:E\to Y$ a topological vector bundle. We apply our observation to $Z=M$ and $V=TM\oplus p^*E\,$. Given embeddings $i_{TM}:TM\to Y\times\reals^k$ and $i_E:E\to Y\times\reals^\ell\,$, we may define an embedding $i_V:V\to E\times\reals^{k+\ell}$ of $V\,$. Thus, we obtain 
\[
(NTM\oplus p^*NE)\times\reals^{k+\ell}\approx\pi_{TM\oplus p^*E}^*((NTM\oplus p^*NE)\otimes\cplxs)\ ,
\]
and a Thom isomorphism 
\[
\vphi_{TM\oplus p^*E}:K^*_c(TM\oplus p^*E)\to K^*_c((NTM\oplus p^*NE)\times\reals^{k+\ell})\ .
\]
Similarly, taking $Z=Y$ and $V=E\,$, and considering the souped-up embedding $E\to Y\times\reals^{2k+\ell}$ given by the embedding of the origin in $\reals^{2k}\,$, we 
obtain a Thom isomorphism $\vphi_E:K^*_c(E)\to K_c^*(NE\times\reals^{2k+\ell})\,$. Since $NTM$ is open in $Y\times\reals^k\,$, we find 
that $(NTM\oplus p^*NE)\times\reals^{k+\ell}$ is open in $NE\times\reals^{2k+\ell}\,$. 

Thus, we obtain a topological family index map (with coefficients in $E$) \`a la Atiyah--Singer \cite{as-ind4} as $\vphi_E^{-1}\circ\vphi_{TM\oplus p^*E}^{\vphantom{-1}}\,$. 
As usual, this object is (as a homomorphism of $K$-groups or an element of $KK$ theory), independent of the choice of embeddings. We shall presently see that this index 
coincides with the Connes--Skandalis map for the previously considered $\mathcal C^{1,0}$ groupoid $\mathcal G=p^*E\times_M(M\times_YM)\,$. 

To this end, we apply the construction detailed in the previous section. We assume now that $p$ is proper, moreover, that the embeddings $i_{TM}$ and $i_E$ are closed, 
and write $i_{TM}=(p,i)$ and $i_E=(\pi_E,j)\,$. Define a homomorphism $h_1:\mathcal G=p^*E\times_M(M\times_YM)\to\reals^{2(k+\ell)}$ by
\[
h_1(x_1,\eta,x_2)=\Parens0{i(x_1)-i(x_2),0,j(p(x_1),\eta),0}\mathfa(x_1,\xi,x_2)\in\mathcal G\ .
\]
This induces a strict homomorphism $h:\mathbb T\mathcal G\to\reals^{2(k+\ell)}\,$,
\[
h(\tau)=
\begin{cases}
  \Parens0{T_xi\xi,0,j(p(x),\eta),0}&\tau=(x,\xi,\eta)\in TM\oplus p^*E\ ,\\
  \Parens0{\eps^{-1}(i(x_1)-i(x_2)),0,j(p(x_1),\eta),0}&\tau=(x_1,\eta,x_2,\eps)\in\mathcal G\times\,\Lopint001\ .
\end{cases}
\]
Then $(r,h,s)$ is injective, and from the closedness of $i,Ti\,$, and $j\,$, it follows it is closed, and hence, proper. Consider the embeddings at $\eps=0$ and 
$\eps=1\,$,
\[
\xymatrix{%
  TM\oplus p^*E=A(\mathcal G)\ar[r]^-{i_0}&\mathbb T\mathcal G&\mathcal G\ar[l]_-{i_1}}
\]
and define $h_0=h\circ i_0\,$. Write $\sigma=\sigma_{\mathbb T\mathcal G,h}\,$, $\sigma_0=\sigma_{A(\mathcal G),h_0}\,$, and $\sigma_1=\sigma_{\mathcal G,h_1}\,$. 
Then we have in $KK$
\[
\sigma\otimes B(i_0)^*=e_{0*}\sigma_0\nd\sigma\otimes B(i_1)^*=e_{1*}\sigma_1\ ,
\]
where $e_0\,$, $e_1$ are the evaluations on $\Cred0{\mathbb T\mathcal G}$ induced by the inclusions $i_0\,$, $i_1\,$. Moreover, as spaces,
\begin{align*}
B_0&=B\Parens0{A(\mathcal G),h_0}=(NTM\oplus p^*NE)\times\reals^{k+\ell}\nda\\
B_1&=B\Parens0{\mathcal G,h_1}=NE\times\reals^{2k+\ell}\ .
\end{align*}
The second homeomorphism is given by 
\begin{gather*}
  B_1=M\times\reals^{2(k+\ell)}/\mathcal G\to NE\times\reals^{2k+\ell}\ ,\\
  [x,u_1,u_2,v_1,v_2]\mapsto\Parens0{p(x),v_1+E_{p(x)},u_1+i(x),u_2,v_2}\ ,
\end{gather*}
and similarly for the first. As noted above, $B_0$ is open in $B_1\,$. Denote the open inclusion by $\iota\,$. Since the proper strict homomorphisms $B(i_0)$ and 
$B(i_1)\circ\iota$ are homotopic, we find $\pmb\iota^*B(i_0)^*=B(i_1)^*$ in $KK\,$, where $\pmb\iota$ is the $*$-morphism $\Ct[_0]0{B_0}\to\Ct[_0]0{B_1}$ 
induced by $\iota\,$.

From the definitions, it follows that $\sigma_0$ is the $KK$ element representing the topological Thom isomorphism $\vphi_{TM\oplus p^*E}\,$. We need to see that 
after applying suitable $KK$ equivalences, the same holds for $\sigma_1\,$. 

To that end, define the strict homomorphism
\[
\pi:\mathcal G\to\mathcal G'=(\reals^n\ltimes\reals^n)\times E:(x_1,\eta,x_2)\mapsto\Parens0{i(x_1)-i(x_2),i(x_2),p(x_1),\eta}\ .
\]
Since $i$ is a closed embedding, $\pi$ is injective and proper. Moreover, it preserves Haar systems if we choose the measure on
$M$ defining the Haar system of $M\times_YM$ to be given by the pullback of the $\dim M^y$-dimensional Hausdorff measure on $\reals^k$ along $i\,$. We denote the 
$*$-morphism $\Cred0{\mathcal G'}\to\Cred0{\mathcal G}$ induced by $\pi$ by the letter $\pmb\pi\,$. 

The homomorphism
\[
h_1':\mathcal G'\to\reals^{2(k+\ell)}:(u,v,e)\mapsto (u,0,j(e),0)
\]
is closed, $(r,h'_1,s)$ is injective, and $h_1=h'_1\circ\pi\,$. We find that 
\[
B(\mathcal G',h_1')=\Parens0{(Y\times\reals^k)\times\reals^{2(k+\ell)}}/\mathcal G'=NE\times\reals^{2k+\ell}=B_1
\]
via
\[
B(\mathcal G',h_1')\to NE\times\reals^{2k+\ell}:[y,u,v_1,v_2,w_1,w_2]\mapsto(y,w_1+E_y,v_1,v_2,w_2)\ ,
\]
and that $B(\pi)=\id_{B_1}\,$.

We have $\Cred0{\mathcal G'}=\Cred0E\otimes\mathbb K$ where $\mathbb K=\Cred0{\reals^k\ltimes\reals^k}\,$. Write $h_1'=j\times\pr_1\,$. There is a canonical 
isomorphism 
\[
\Cred0{\reals^k\ltimes\reals^k}\rtimes_{\alpha_{\pr_1}}\reals^{2k}\to\Ct[_0]0{\reals^{2k}}\rtimes\reals^{2k}=\mathbb K
\]
where the latter action is regular (i.e.~dual to the trivial action). Because the Thom element in $KK(\Ct[_0]0{\reals^{2k}},\Ct[_0]0{\reals^{2k}}\rtimes\reals^{2k})$ is 
the Bott element \cite[Example 19.3.4 (ii)]{blackadar-ktheory}, composition of this isomorphism with $t_{\alpha_{\pr_1}}$ is the identity in $KK(\mathbb K,\mathbb K)=\ints\,$. 

Thus, the Thom element $t_1'=t_{\alpha_{h_1'}}$ defining $\sigma_1'=\sigma_{\mathcal G',h_1'}$ satisfies $\pi^*t_1=t_1'$ where 
$t_1=t_{\alpha_{h_1}}\,$, and on the other hand, $\sigma_1'$ induces on $K$-theory the topological Thom isomorphism $\vphi_E\,$. On the level of $KK$ theory, 
\[
\tau=\sigma_0\otimes\sigma_1^{-1}=\sigma_0\otimes\pmb\pi_*(\sigma_1')^{-1}=\pmb\pi_*\Parens0{\vphi_{TM\oplus p^*E}^{\vphantom{-1}}\otimes\vphi_E^{-1}}
\]
Thus, we have proved the following theorem.

\begin{thm}\label{cs-top-expr}
  Let $E\to Y$ be a topological vector bundle and $p:M\to Y$ a class $\mathcal C^{1,0}$ manifold \scth $p$ is proper. Then the Connes Skandalis map $\tau$ associated
  to the tangent groupoid of $\mathcal G=p^*E\times_M(M\times_YM)$ is 
  \[
  \tau=\pmb\pi_*\Parens0{\vphi_{TM\oplus p^*E}^{\vphantom{-1}}\otimes\vphi_E^{-1}}
  \]
  where $\mathbf{\pi}:\Cred0E\otimes\mathbb K\to\Cred0{\mathcal G}$ is a $KK$ equivalence and $\vphi_{TM\oplus p^*E}^{\vphantom{-1}}\otimes\vphi_E^{-1}$ is the
Atiyah--Singer topological family index for $p^*E\oplus TM\,$. 
\end{thm}

\section{Proof of the Wiener--Hopf Index Formula}

  In this section, we shall prove the index formula for the Wiener--Hopf C$^*$-algebra explained in the introduction. Recall from the introduction and from 
  \cite{alldridge-johansen-wh1} that on the level of $KK$ theory, this was an expression of the $KK^1$ element $\partial_j$ representing the extension
  \[
  \xymatrix{%
    0\ar[r]&\Cred0{\mathcal G|U}\ar[r]&\Cred0{\mathcal G}\ar[r]&\Cred0{\mathcal G|F}\ar[r]&0}
  \]
  where $\mathcal G=\mathcal W_\Omega|(U_{j+1}\setminus U_{j-1})$ is the restriction of the Wiener--Hopf groupoid to the locally closed invariant subset
  $U_{j+1}\setminus U_{j-1}=Y_j\cup Y_{j-1}\,$, and the closed invariant subset $F=Y_j=U_j\setminus U_{j-1}$ is the stratum of the Wiener--Hopf 
  compactification $\clos\Omega$ corresponding to the orbit type of $P_j\,$, and consists of all points of $\clos\Omega$ lying above $P_j\,$, the set of 
  $n_{d-j}$-dimensional faces.

  Our proof proceeds in three steps:{\def\theenumi{\arabic{enumi}}
  \begin{enumerate}
  \item Construct a class $\mathcal C^{1,0}$ groupoid $\mathcal D_j\,$, which is of the form 
    \[
    \mathcal D_j=p^*E\times_M(M\times_YM)
    \]
    for some $\mathcal C^{1,0}$ manifold $M\to Y\subset P_{j-1}$ and some vector bundle $E\to Y\,$.
  \item Construct a proper strict homomorphism $\mathbb W\mathcal D_j\to\mathcal W_\Omega|(U_{j+1}\setminus U_{j-1})\,$, where $\mathbb W\mathcal D_j$ is the
    `suspended tangent groupoid' constructed in section \ref{suspcone-tangent}. By naturality, $\partial_j$ is expressed in terms of the standard extension belonging
    to $\mathbb W\mathcal D_j\,$.
  \item Now we are in a position to apply the $KK$ theoretical yoga explained above. By the results of \ref{susp-conneskand}, we relate the extension induced by
    $\mathbb W\mathcal D_j$ to the Connes--Skandalis map pertaining to the tangent groupoid of $\mathcal D_j\,$. Because of the particular form of the latter groupoid,
    the computations from Section \ref{topexpr-connesskand} furnish the $KK$ theoretical index formula.
  \end{enumerate}}

\subsection{Embedding the Fibres of $\mathcal P_j$}

  Returning to the study of the Wiener--Hopf groupoid $\mathcal W_\Omega\,$, we shall take the first of the three steps mentioned in this section's introduction, comprising
  the proof of the Wiener--Hopf index formula. 

  Recall to that end the notation and the notions from \cite{alldridge-johansen-wh1}. In particular, $\Omega$ is a pointed, solid, 
  closed convex cone in the finite dimensional Euclidean vector space $X$ ($=\reals^n$). Let $\{0=n_0<n_1<\dotsm<n_d=n\}$ be the set of dimensions of the faces of 
  the dual cone $\Omega^*\,$, ordered increasingly (compare the introduction). Set 
  \[
  P_j=\Set1{F\subset\Omega^*}{\dim F=n_{d-j}\,,\,F\text{ face of }\Omega^*}\ .
  \]
  We shall assume that the cone $\Omega^*$ is \emph{facially compact}, i.e.~all the spaces $P_j$ of are compact in the topology induced from the space $\closed0X$ of all 
  closed subsets of $X\,$, endowed with the Fell topology. Consider
  \[
  \mathcal P_j=\Set1{(E,F)\times P_{j-1}\times P_j}{E\supset F}\ ,
  \]
  which is a compact subspace of $P_{j-1}\times P_j\,$. We have projections $\xi:\mathcal P_j\to P_{j-1}$ and $\eta:\mathcal P_j\to P_j\,$. 

  We shall prove that under suitable assumptions, $\mathcal P_j$ is in a natural way a $\mathcal C^{1,0}$ manifold with respect to the projection $\xi$ onto the base 
  $\xi(\mathcal P_j)\subset P_{j-1}\,$, and moreover, we will give a specific Euclidean embedding. It is important to stress that the mere existence of a $\mathcal C^{1,0}$ 
  structure is not sufficient for our needs: to ultimately evaluate the Atiyah--Singer family index map associated to the groupoid to be constructed, a concrete Euclidean 
  embedding has to be at hand. 

  The $\mathcal C^{1,0}$ manifold $\mathcal P_j$ forms the basic building block for our sought-for $\mathcal C^{1,0}$ groupoid $\mathcal D_j\,$. Indeed, as groupoids,  
  $\mathcal D_j=\xi^*\Sigma_{j-1}\times_{P_{j-1}}(\mathcal P_j\times_{P_{j-1}}\mathcal P_j)\,$. Here, $\Sigma_{j-1}$ is the vector bundle which appears 
  in the construction of the composition series of the Wiener--Hopf C$^*$-algebra (compare the indroduction).

  To introduce a $\mathcal C^{1,0}$ structure on $\mathcal P_j\,$, we first construct suitable Euclidean embeddings of the fibres of $\xi\,$, and show that these admit
  tangent spaces of fixed dimension at every point. 

  To that end, we adopt the following notation: For $A\subset X\,$, let $\Span0A$ denote the linear span of $A\,$, $A^\perp$ the orthogonal complement, $A^*$ the dual cone, 
  and $A^\circledast=A^*\cap\Span0A$ the relative dual cone. Then, for $(E,F)\in\mathcal P_j\,$, define 
  \[
  E_{1/2}(F)=F^\perp\cap(F^\perp\cap E^\circledast)^\perp\ .
  \]
  We shall be using some convex geometry in the sequel, and refer the reader to \cite[Chapter I]{hilgert-hofmann-lawson} for a thorough treatment. 

  \begin{rem}
    Although the notation $E_{1/2}(F)$ may seem somewhat arbitrary, we have chosen it to stress the analogy to the case of symmetric cones. We briefly explain its  
    meaning in this case, and refer to \cite{faraut_koranyi} for details on this topic. 

  Indeed, assume that $\Omega=\Omega^*$ is a symmetric cone in the Euclidean vector space $X\,$, i.e.~a self-dual cone whose interior admits a transitive action 
  of a subgroup of $\GL(X)\,$. Then, $X$ is in an essentially unique way a Euclidean Jordan algebra. If $c\in X$ is an idempotent, then we have the orthogonal 
  `Peirce decomposition' 
  \[
  X=X_0(c)\oplus X_{1/2}(c)\oplus X_1(c)\mathtxt{where}X_\lambda(c)=\ker(L(c)-\lambda)
  \]
  are the eigenspaces for the action of $c$ on $X$ by left multiplication, called `Peirce spaces', and $\dim X_1(c)$ is called the rank of $c\,$. 

  The faces of $\Omega$ are of the form $E=\Omega\cap X_0(e)\,$, for idempotents $e\in X\,$. Now assume that $E\supsetneq F$ are faces \scth $F$ has 
  minimal codimension in $E\,$. Then $E=\Omega\cap X_0(e)\,$, $F=\Omega\cap X_0(c)$ where $e,c$ are idempotents with $c-e$ an idempotent of rank
  one. The dual face of $F$ in the self-dual cone $E$ is the extreme ray $E\cap X_1(c-e)\,$. The Euclidean Jordan algebra $\Span0E=X_0(e)$ has the Peirce decomposition
  \[
  X_0(e)=X_0(c)\oplus\Parens0{X_0(e)\cap X_{1/2}(c-e)}\oplus X_1(c-e)
  \]
  w.r.t.~$c-e\in X_0(e)\,$. Here, $X_0(c)=\Span0F$ is the linear span of the face $F\,$, $X_1(c-e)$ is the line spanned by dual cone of $F$ in $E\,$, and 
  $X_0(e)\cap X_{1/2}(c-e)$ is the intersection of the orthogonal complements of the two former spaces. 

  Moreover, the set of all proper faces $F$ of $E$ of minimal
  codimension corresponds exactly to the set $S$ of rank one idempotents of the Euclidean Jordan algebra $X_0(e)\,$. The set $S$ is the `Shilov boundary' of
  the bounded domain whose Siegel realisation is $X_0(e)+i\Omega\cap X_0(e)\,$. It is a compact submanifold of $X_0(e)\,$, and the tangent space at an idempotent 
  $f\in S$ is precisely $X_0(e)\cap X_{1/2}(f)\,$.
\end{rem}

  Returning to the general case of a no longer necessarily symmetric cone $\Omega\,$, we shall see that under mild conditions, the Peirce decomposition that we have explained 
  for the case of symmetric cones has a counterpart for any pair $(E,F)\in\mathcal P_j\,$, the fibre of $\xi$ over $E$ corresponds exactly to a compact set of generators 
  of extreme rays in $E^\circledast\,$, and the space $E_{1/2}(F)$ can be interpreted as the tangent space of this set at the point corresponding to $F\,$. First, let us 
  explain what kind of geometric conditions have to be imposed on $\Omega^*\,$.

\begin{defn}
  A face $E$ of a cone $C$ will be called \emph{modular} if it contains a face $F$ which is maximal dimensional relative $C\,$, i.e.
  \[
  \dim F=\max\Set1{\dim G}{G\text{ a face of }C\,,\,\dim G<\dim E}\ .
  \]
\end{defn}

\begin{rem}
Applied to the cone $C=\Omega^*\,$: For $j\sge1\,$, $\xi(\mathcal P_j)$ is the collection of a modular faces of dimension $n_{d-j+1}\,$. 
\end{rem}

\begin{lem}\label{fibrewise-tangents}
  Let $E$ be a pointed cone, and $F\subsetneq E$ a maximal face. Then $F$ is exposed in $E$ (i.e.~equals the intersection of $E$ with a supporting hyperplane), and the relative 
  dual face $F^\perp\cap E^\circledast\neq0$ is an extreme ray of $E^\circledast$ if and only if it contains an exposed ray of $E^\circledast\,$. In this case, we 
  denote the normalised generator of $F^\perp\cap E^\circledast$ by $e_F=e_F(E)\,$.
\end{lem}

\begin{pf}
  There exists a face $F\subsetneq G\subset E$ \scth $F$ is exposed in $G\,$. Since $F$ is maximal, $G=E\,$, and $F$ is exposed in $E\,$. The dual face
  $F^\perp\cap E^\circledast$ is an exposed face of $E^\circledast\,$. So, if it is an extreme ray, then it is exposed.

  Conversely, let $G\subset F^\perp\cap E^\circledast$ be an exposed ray of $E^\circledast\,$. Then $F\subset E\cap G^\perp\,$. 
  Of course, $E\cap G^\perp\neq E\,$, since $\dim G^\perp<\dim E\,$. The maximality of $F$ entails $F=E\cap G^\perp\,$. Since $G$ is exposed in $E^\circledast\,$,
  it follows that $G=F^\perp\cap E^\circledast\,$, by \cite[Proposition I.2.5]{hilgert-hofmann-lawson}, and hence, an extreme ray.
\end{pf}

\begin{defn}
  Let $E$ be a modular face of a cone $C\,$. We shall say that $E$ is \emph{smooth} if the relative interior of every face $F\subsetneq E\,$, maximal dimensional 
  relative $C\,$, consists of regular points of $E\,$. Here, a point of $E$ is called regular if it admits a unique supporting hyperplane, 
  cf.~\cite[Definition I.2.24]{hilgert-hofmann-lawson}. 

  We shall say that a cone is \emph{locally smooth} if all of its modular faces are smooth. 
\end{defn}

\begin{rem}
  Let $E\subset C$ be a modular face. By Lemma \ref{fibrewise-tangents} and \cite[Proposition I.2.25]{hilgert-hofmann-lawson}, $E$ is smooth if and only if for all faces 
  $F\subset E$ maximal dimensional relative $C\,$, all the extreme rays of the relative dual face $F^\perp\cap E^\circledast$ are exposed in $E^\circledast\,$. In particular, 
  if $E^\circledast$ is facially exposed (i.e.~all faces are exposed), then $E$ is smooth. (However, the condition that all extreme rays of a cone be exposed does not imply 
  that all faces are exposed.)
\end{rem}

\begin{cor}\label{cone-embedding}
  Let the cone $\Omega$ be such that $\Omega^*$ is locally smooth. Then for each $1\sle j\sle d$ there is a well-defined injection
  \[
  \mathcal P_j\to P_{j-1}\times X:(E,F)\mapsto\Parens0{E,e_F(E)}\ .
  \]
\end{cor}

We conclude this subsection by a brief digression to gauge the generality of the condition of local smoothness.

\begin{lem}\label{locsmooth-cones}
The following classes of cones are locally smooth and have locally smooth dual cones:
\begin{enumerate}\romannum
\item Polyhedral cones,
\item Lorentz cones, and
\item homogeneous cones, in particular, symmetric cones.
\end{enumerate}
\end{lem}

\begin{pf}
(i). The dual of a polyhedral cone is polyhedral \cite[Corollary I.4.4]{hilgert-hofmann-lawson}. Any face of a polyhedral cone is polyhedral, and polyhedral cones are
  facially exposed.

(ii). The dual of Lorentz cone is a Lorentz cone, and its non-zero proper faces are all exposed extreme rays, cf.~\cite[Proposition I.4.11]{hilgert-hofmann-lawson}.

(iii). The dual of a homogeneous cone is itself homogeneous \cite[Satz~4.3]{dorfmeister-koecher-regcone}. Homogeneous cones are facially exposed 
  \cite[Theorem 3.6]{truong-tuncel}. So, it remains to see that all faces of homogeneous cones are homogeneous. To see this, we recall Rothaus's \cite{rothaus-constrhomcones} 
  inductive construction of homogeneous cones. If $U,V$ are finite-dimensional vector spaces, $K\subset V$ a closed convex cone, and $B:U\times U\to V$ 
  is a bilinear map, then we say that $B$ is $K$-\emph{positive} if $B(u)=B(u,u)\in K\setminus0$ \fa $u\in U\setminus0\,$. Given such data, the \emph{Siegel cone}
  \[
  C(K,B)=\Set1{(u,v,t)\in U\times K\times\reals_{\sge0}}{tv-B(u)\in K}
  \]
  is a closed convex cone in $U\times V\times\reals\,$. If $K$ is homogeneous, then the bilinear map $B$ is called \emph{homogeneous} if for some subset $G\subset GL(V)$
  acting transitively on $K^\circ\,$, and \fa $g\in G\,$, there exist elements $g_U\in\GL(U)$ \scth
  \[
  g B(u,u')=B\Parens0{g_U(u),g_U(u')}\mathfa u,u'\in U\ .
  \]
  If $K$ and $B$ are homogeneous, then $C(K,B)$ is homogeneous. Conversely,
  any homogeneous cone $C$ can be obtained by this procedure from a homogeneous cone $K$ of dimension strictly less than $\dim C\,$. Thus, all homogeneous cones can be 
  constructed inductively from the real half-line, cf.~\cite{vinberg-convexhomogeneous,rothaus-constrhomcones,truong-tuncel}. All faces of the half-line are homogeneous. So, we need to see
  that this property remains stable under the inductive step of the above construction.

  By \cite[Theorem 3.2]{truong-tuncel}, the set $\ext C(K,B)$ of generators of extreme rays for a homogeneous Siegel cone is given as follows:
  \[
  \ext C(K,B)=\Set1{(u,v,t)\in U\times K\times\reals_{\sge0}}{tv=B(u)\,,\,t=0\,\imp\,v\in\ext K}\ .
  \]
  Now suppose that $E\subset C(K,B)$ is a face, and let $F=(0\times K\times 0)\cap E\,$, which defines a face of $K\,$. Then, by \cite[proof of Theorem 3.6]{truong-tuncel}, either
  $E$ contains only extreme generators of the form $(0,v,0)\,$, $v\in\ext K\,$, or we have the equivalence $B(u)\in F\,\equiva\,(u,B(u),1)\in E\,$. Any cone is the positive linear
  span of its extreme generators, by \cite[Theorem I.3.16]{hilgert-hofmann-lawson}. Thus, in the former case, $E=F\,$, in which case we are done by our inductive hypothesis. In the
  latter case,
  \begin{align*}
  \ext E&=\Set1{(u,v,t)\in\ext C(K,B)}{B(u)\in F}\\
  &=\Set1{(u,v,t)\in U\times F\times\reals_{\sge0}}{tv=B(u)\,,\,t=0\,\imp\,v\in\ext F}\ .
  \end{align*}
  Then, define $U_F=\Set0{u\in U}{B(u)\in F}\,$. This set is a linear subspace of $U\,$. Indeed, if $u,v\in U_F\,$, then
  \[
  \tfrac12\cdot(B(u+v)+B(u-v))=B(u)+B(v)\in F+F=F\ .
  \]
  Since $B(u+v),B(u-v)\in K$ and $F\subset K$ is a face, it follows that $B(u\pm v)\in F\,$, so $u\pm v\in U_F\,$. Since $U_F$ is invariant under positive scalar multiples,
  it is indeed a linear subspace.

  But then $B_F=B|(U_F\times U_F)$ is $F$-positive, and $C(F,B_F)$ makes sense. In fact, $\ext E$ is the set of extreme generators of $C(F,B_F)$ by our previous calculations,
  as soon as we have established that $B_F$ is homogeneous. In that case, it will follow that $E=C(F,B_F)\,$, both being the positive linear spans of their extreme generators,
  thereby establishing the homogeneity of $E\,$.

  So, let us check that $B_F$ is homogeneous. By our inductive assumption, $F$ is homogeneous, and we may choose a subset $G\subset\GL(\Span0F)$ acting transitively on
  $F^\circ\,$. Since $B$ is homogeneous, to $g\in G\,$, there exists $g_U\in\GL(U)$ \scth $gB(u,u')=B(g_U(u),g_U(u'))\,$. If $u\in U_F\,$, then
  \[
  B\Parens0{g_U(u)}=g\Parens0{B(u)}\in g(F)=F\ ,
  \]
  so $g_U$ leaves $U_F$ invariant, and $B_F$ is homogeneous. This proves our claim.
\end{pf}

\begin{figure}[h]
  \centering
  \begin{center}
    \includegraphics[bb=0 0 274 126,scale=0.5]{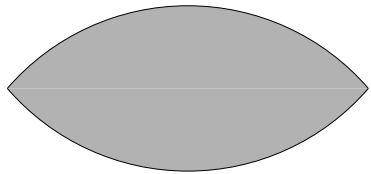}
  \end{center}
  \caption{The section of a three-dimensional cone which is not locally smooth but whose dual cone is.}
  \label{fig:almond}
\end{figure}

\begin{rem}
  Although $\Omega^*$ is certainly locally smooth in the most interesting cases, let us note in passing the curious asymmetry of this condition. Indeed, consider the 
  three-dimensional, facially exposed cone $C$ which has the `almond shaped' section illustrated in Fig.~\ref{fig:almond}. $C$ itself not locally smooth, since 
  it admits several supporting hyperplanes at the tips of the almond. The dual cone $C^*\,$, however, is locally smooth since its non-zero and proper modular faces are 
  polyhedral.
\end{rem}

\subsection{A $\mathcal C^{1,0}$ Structure on $\mathcal P_j$}

  Having verified that local smoothness of the dual cones obtains for large classes of cones, we shall in the sequel always \emph{assume that $\Omega^*$ is locally smooth}. In
  particular, Lemma \ref{fibrewise-tangents} allows for the definition of extremal generators $e_F(E)=e_F\in E^\circledast$ for any $(E,F)\in\mathcal P_j\,$, and by 
  Corollary \ref{cone-embedding}, this defines a set-theoretical embedding of the fibres of $\xi\,$. 

  We shall show that the embedding is in fact fibrewise bi-Lipschitz, and that the image is a closed $\mathcal C^{1,0}$ submanifold of $\mathcal P_j\times X\,$. To accomplish 
  this, we apply a simple but effective idea used by Pugh in his dynamical systems proof of the Cairns--Whitehead Theorem \cite[Proof of Proposition 3.3]{pugh-cwtheorem}. 
  Namely, it is a basic fact from analysis that if the derivative of a Lipschitz map $g$ (which exists almost everywhere due to Rademacher's Theorem) can be extended 
  continuously, then $g$ is $\mathcal C^1\,$. Thus, if $\xi:\mathcal P_j\to P_{j-1}$ is given a fibrewise Lipschitz manifold structure, it is sufficient to show 
  that the family of the tangent spaces which exist almost everywhere in every Lipschitz chart, extends to a continuous map to an appropriate 
  Grassmannian. It follows that the fibres have regular local $\mathcal C^1$ parametrisations. 

  The main part of the proof will in fact consist in establishing the Lipschitz structure. 

\begin{prop}\label{fibrewise-tangents2}
  The map
  \[
  \mathbf e:\mathcal P_j\to X:(E,F)\mapsto e_F(E)
  \]
  is a closed embedding and bi-Lipschitz on every fibre of $\xi:\mathcal P_j\to P_{j-1}\,$, locally uniformly with respect to the fibres. Here, the metric on
  $\mathcal P_j$ is the box metric for a metric on the space $P=\bigcup_{i=0}^dP_i$ inducing the Fell topology on this set.
\end{prop}

In the proposition's proof, we shall need a metric which defines the Fell topology on the set $\mathcal C$ of all closed convex cones in $X\,$. For all closed subsets $A,B\subset X\,$, define the excess functional $e(A,B)=\sup\Set0{\dist(a,B)}{a\in A}\,$, and set 
\[
h(A,B)=\max\Parens0{e(A\cap\ball,B),e(B\cap\ball,A)}\mathfa A,B\in\mathcal C
\]
where $\ball$ is the unit ball in $X\,$. If $a\in A\cap\ball\,$, then the projection $b=\pi_B(a)$ onto $B$ has $\Norm0b\sle\Norm0a\sle1\,$, so $\dist(a,B)=\dist(a,B\cap\ball)\,$.
Hence, $h=H\,$, where 
\[
H(A,B)=\max\Parens0{e(A\cap\ball,B\cap\ball),e(B\cap\ball,A\cap\ball)}\mathfa A,B\in\mathcal C\ .
\]
Since $H(A,B)$ is the Hausdorff distance of $A\cap\ball$ and $B\cap\ball\,$, $h$ defines a metric. By \cite[Exercise 5.1.10, Lemma 7.2.6]{beer-topologies-book}, the topology induced 
on $\mathcal C$ by the distance function $h$ coincides with the subspace topology from $\closed0X\,$. Moreover, the Walkup--Wets Isometry Theorem 
\cite[Theorem 7.2.9]{beer-topologies-book} shows that the map $C\mapsto C^*$ is an isometry for $h\,$.

{\def\Elproofname{\protect{PROOF of Proposition \ref{fibrewise-tangents2}.}}
\begin{pf}
  Observe that the map $\mathbf e$ has closed graph. In fact, given sequences $(E_k,F_k)\to(E,F)$ and $e_{F_k}(E_k)\to e\,$, then by \cite[Proposition 2.2.8]{alldridge-johansen-wh1}
  and continuity of polarity \cite[Corollary 7.2.12]{beer-topologies-book}, we have $\Span0{E_k}\to\Span0E\,$, $F_k^\perp\to F^\perp\,$, and $E_k^*\to E^*\,$. From the
  definition of Painlev\'e--Kuratowski convergence (cf.~\cite[Section 2.1]{alldridge-johansen-wh1}), we conclude 
  \[
  e\in\Span0E\cap E^*\cap F^\perp=\reals_{\sge0}\cdot e_F(E)\ .
  \]
  Since $e$ is a unit vector, we have $e=e_F(E)\,$. Thus, $\mathbf e$ indeed has closed graph. Because $\mathcal P_j$ and the unit sphere of $X$ are compact, it follows that
  $\mathbf e$ is continuous and closed.

  Now, to the bi-Lipschitz continuity on every fibre. Fix $(E,F)\in\mathcal P_j\,$. Because $\mathbf e$ is continuous and
  $e_F$ is a unit vector, there exists an open neighbourhood $U_{E,F}\subset\mathcal P_j$ of $(E,F)\,$, \scth
  \[
  \Rscp0{e_{H_1}(G_1)}{e_{H_2}(G_2)}>0\mathtxt{for every $(G_j,H_j)\in U_{E,F}\,$, $j=1,2\,$.}
  \]
  Let $(G,H_j)\in U_{E,F}\,$, $j=1,2\,$. Writing $e_j=e_{H_j}\,$, we obtain
  \[
    h(H_1,H_2)=h(H_1^*\cap\Span0G,H_2^*\cap\Span0G)=\Norm0{e_1-\Rscp0{e_1}{e_2}\cdot e_2}=\sqrt{1-\Rscp0{e_1}{e_2}^2}\ ,
  \]
  because $C\mapsto C^*$ is an isometry and we may consider $h$ relative to $\Span0G\,$. Since $\Rscp0{e_1}{e_2}\in[0,1]\,$, we have  
  \[
  \sqrt2\cdot\sqrt{1-\Rscp0{e_1}{e_2}^2}\sge\Norm0{e_1-e_2}=\sqrt{2-2\Rscp0{e_1}{e_2}}\sge\sqrt{1-\Rscp0{e_1}{e_2}^2}\ .
  \]
  Thus, the map $\mathbf e$ is bi-Lipschitz when restricted to $U_{E,F}\cap\xi^{-1}(G)\,$, with Lipschitz constants
  independent of $G\in\xi(U_{E,F})\,$. By compactness, we obtain global Lipschitz conditions.
\end{pf}}

\begin{thm}\label{c10-field}
The image of $\mathcal P_j$ under $(\xi,\mathbf e)$ is a closed $\mathcal C^{1,0}$ submanifold of $\xi(\mathcal P_j)\times X$ (which is a $\mathcal C^{1,0}$ manifold over $\xi(\mathcal P_j)$ with respect to the first projection).  This defines a $\mathcal C^{1,0}$ manifold structure on $\mathcal P_j$ over $\xi(\mathcal P_j)\subset P_{j-1}$ \scth the fibrewise tangent space at $(E,F)\in\mathcal P_j$ is $E_{1/2}(F)\,$.
\end{thm}

The theorem's proof follows the basic strategy outlined at the beginning of the subsection and decomposes into a sequence of lemmata.

\begin{lem}\label{shilov-tangent}
Let $(E,F),(E,F^i_k)\in\mathcal P_j\,$, $F_k^i\to F\,$, $i=0,1\,$, and $\eps_k\to0+\,$. Whenever the limit
$v=\limk_k\eps_k^{-1}\cdot(e_{F_k^0}-e_{F_k^1})$ exists, it lies in $E_{1/2}(F)\,$.
\end{lem}

\begin{pf}
  First, we have $v\perp F\,$. Indeed, write $e_k^i=e_{F_k^i}\,$. If $f\in F\,$, then there exist $f_k^j\in F_k^i$ such that $f_k^i\to f\,$. In particular,
  \[
  \Rscp0vf=\limk_k\eps_k^{-1}\cdot\Rscp0{e_k^0-e_k^1}{f_k^i}=(-1)^{1-i}\cdot\limk_k\eps_k^{-1}\cdot\Rscp0{e_k^{1-i}}{f_k^i}\sge0\ ,
  \]
  since $e_k^{1-i}\in E^\circledast\subset\Omega$ and $f_k^i\in F_k^i\subset E\,$, so $\Rscp0vf=0\,$. This shows that $v\perp F\,$, as desired. Moreover, 
  it is clear that $v\in\Span0E\,$.

  Remains to prove that $v\perp e_F\,$. This is seen similarly, namely
  \[
  \Rscp0v{e_F}=\limk_k\eps_k^{-1}\cdot\Rscp0{e_k^0-e_k^1}{e_k^i}=(-1)^i\cdot\limk_k\eps_k^{-1}\cdot\Parens0{1-\Rscp0{e_k^0}{e_k^1}}\ ,
  \]
  and thus $\Rscp0v{e_F}=-\Rscp0v{e_F}=0\,$. This proves our assertion.
\end{pf}

\begin{lem}\label{transverse}
Given $(E,F)\in\mathcal P_j\,$, the map $p=p_{E_{1/2}(F)}$ is injective in a neighbourhood $U$ of $e_F=e_F(E)$ in $\mathbf e(\xi^{-1}(E))\,$. 
\end{lem}

\begin{pf}
  To establish the injectivity of $p\,$, assume, seeking a contradiction, that for each neighbourhood $U\subset \mathbf e(\xi^{-1}(E))$ of $e_F\,$, there exist $y_1,y_2\in U\,$, $y_1\neq y_2\,$, \scth $y_1-y_2\in E_{1/2}(F)^\perp\,$. Then there are sequences $y_j^k\in \mathbf e(\xi^{-1}(E))\,$, $j=1,2\,$, \scth 
  \[
  0<\Norm0{y_1^k-y_2^k}\sle\tfrac1k\ ,\ y_1^k-y_2^k\perp E_{1/2}(F)\ . 
  \]
  Passing to a subsequence, the limit $v=\limk_k\Norm0{y_1^k-y_2^k}^{-1}(y_1^k-y_2^k)$ exists. Then $v\in E_{1/2}(F)^\perp\,$, and by Lemma \ref{shilov-tangent}, 
  $v\in E_{1/2}(F)\,$. This is a contradiction, since $\Norm0v=1\,$.  
\end{pf}

\begin{lem}\label{mapspropballimage}
  Let $K\subset\reals^n$ be a compact subset, $B\subset\reals^m$ a closed ball centred at the origin. Assume that $f,g:K\to B$ are 
  continuous, and that \fa $x\in K\,$, we have $f(x)=\lambda\cdot g(x)$ for some $\lambda>0$ (depending on $x$). If $g$ has convex fibres and $g(K)=B\,$, then $f(K)$ 
  contains $\delta\cdot B$ for some $\delta>0\,$. 
\end{lem}

\begin{pf}
  In fact, there is no restriction in assuming that $B$ is the unit ball. The map $g$ is continuous and open, and its fibres are compact convex sets. Thus, by Michael's Selection 
  Theorem \cite[Theorem 3.1']{michael-sel1}, $g$ has a continuous section, $s\,$, say. 

  Let $x\in B\,$, $\Norm0x=1\,$. Then $(f\circ s)(0)=0$ and 
  $(f\circ s)([0,x])\subset[0,x]\,$. Since $f\circ s$ has some non-zero value on $[0,x]\,$, because $g$ does, we deduce that 
  \[
  \delta_x=\mathop{\mathrm{diam}} (f\circ s)([0,x])>0\ ,
  \]
  and $[0,\delta_x\cdot x]=(f\circ s)([0,x])\subset f(K)\,$. It is clear that $x\mapsto\delta_x$ is lower semi-continuous and thus assumes its infimum. Thus, there 
  exists $\delta>0$ such that $[0,\delta x]\subset f(K)$ \fa $x\in B\,$, $\Norm0x=1\,$. Hence, $\delta\cdot B\subset f(K)\,$. 
\end{pf}

\begin{rem}\label{normal-projection}
  Let $C\subset\reals^n$ be a closed convex set, $U_0\subset\reals^n$ a linear subspace, $u\in U_0^\perp\,$, $U=U_0+u\,$, and $x\in\partial(C\cap U)\,$. If $y\in X$ is 
  an inner normal to $C$ in $x\,$, i.e.~$\Rscp0xy=\inf_{z\in C}\Rscp0zy\,$, then $p_{U_0}(y)$ is an inner normal to $p_{U_0}(C\cap U)=C\cap U-u$ in $p_{U_0}(x)=x-u\,$. 
  In fact, by the Projection Theorem, $\Rscp0zy=\Rscp0{z-u}{p_{U_0}(y)}\,$.
\end{rem}

\begin{lem}\label{proj-open}
  Let $(E,F)\in\mathcal P_j\,$. The map $p=p_{E_{1/2}(F)}:\mathbf e(\xi^{-1}(E))\to E_{1/2}(F)$ is open in a neighbourhood of $e_F=e_F(E)\,$. 
\end{lem}

\begin{pf}
  The map $p$ is certainly open onto its image, so it suffices to show that the image contains some ball. \def\pmax{\partial_{\mathrm{max}}}

  Denote by $\pmax E$ the set of all those boundary points of $E$ which lie in the relative interior of face of $E$ which is maximal dimensional relative $\Omega^*\,$. 
  Then $\pmax E$ is an open subset of $\partial E\,$, and consists of regular points of $E\,$. In particular, the Gau\ss{} map $\nu:\pmax E\to X$ which 
  to $x\in\pmax E$ associates the unique inner unit normal to $E$ in $x\,$, is well-defined. The fibre $\xi^{-1}(E)$ is precisely the quotient of $\pmax E$ 
  by the equivalence relation induced by $\nu\,$, and $\nu$ drops to the map $\xi^{-1}(E)\to \mathbf e(\xi^{-1}(E)):G\mapsto e_G(E)\,$.

  Hence, it is sufficient to see that the image of the restriction of $p\circ\nu$ to some neighbourhood of a point in the relative interior of $F$ contains some ball. Fix 
  $x\in F^\circ$ and let $e=e_F(E)\,$. Then the plane $\Span0{x,e}$ intersects $E^{\circledast\circ}\cap E^\circ$ in some point $h^*$ \scth $\Rscp0{h^*}x=1\,$. 

  Let $S_0=h^{*\perp}$ and $S=S_0+h$ where $h=\Norm0{h^*}^{-2}\cdot h^*\,$. Then $E\cap S$ is a compact convex base of $E\,$. Since 
  $h\in\Span0{x,e}\subset E_{1/2}(F)^\perp\,$, we have that $E_{1/2}(F)\subset S_0\,$. 

  Since $h\in E^{\circledast\circ}$ and $E^{\circledast\circ}\cap F^\perp=\vvoid\,$, $p_{S_0}|F^\perp$ is a linear isomorphism of $F^\perp$ onto $T_0=p_{S_0}(F^\perp)\,$. 
  Thus, $E_{1/2}(F)\subset T_0$ is a hyperplane, and $p$ factors through $p_{T_0}\,$. Let $T=T_0+\tilde x\subset S_0\,$. Then $\tilde x=p_{S_0}(x)=x-h$ is an exposed 
  point of the compact convex set $C=(E-h)\cap T\,$. Indeed, $\tilde e=p_{S_0}(e)$ is an inner normal at $\tilde x$ by Remark \ref{normal-projection}, so 
  $\tilde x+E_{1/2}(F)=\tilde x+S_0\cap\tilde e^\perp$ is a supporting hyperplane for $C\,$, and $C\cap(\tilde x+E_{1/2}(F))=\tilde x\,$. 

  Let $\pmax C=T\cap(\pmax E-h)\,$. Then the Gau\ss{} map $\nu_0:\pmax C\to T_0$ which to any point $y$ associates the unique 
  unit inner normal to $C$ in $y\,$, is well-defined. We claim that for any $y\in\pmax C$ sufficiently close to $\tilde x\,$, there exists $\lambda>0$ \scth 
  $(p\circ\nu)(y+h)=\lambda\cdot(p\circ\nu_0)(y)\,$. If this is the case and the image of $p\circ\nu_0$ on some compact subset of $\pmax C$ contains some ball, 
  then by Lemma \ref{mapspropballimage}, the same is true for the image of $p\circ\nu\,$. Indeed, the fibres of $\nu_0$ are faces of $C\,$, and $p$ is locally injective 
  by Lemma \ref{transverse}. 

  By Remark \ref{normal-projection}, \fa $y\in\pmax C\,$, we have $p_{T_0}(\nu(y+h))=\lambda\cdot\nu(y)$ for some $\lambda\sge0\,$, and we may apply $p$ on both sides of 
  the equation. If $y=\tilde x=x-h\,$, then $\nu(x)=e\in F^\perp\,$, and $p_{T_0}|F^\perp$ is an isomorphism, so $p_{T_0}(\nu(x))\neq0\,$. Thus, $\lambda>0$ if $y$ 
  is sufficiently close to $\tilde x\,$. 

  We are now reduced to proving that $(p\circ\nu_0)(K)$ contains some ball for some compact neighbourhood $K\subset\pmax C$ of $\tilde x\,$. This follows similarly as in 
  \cite[Proof of Lemma 5.2]{fedotov}. Since we are not ware of an English translation of Fedotov's paper, we detail the argument. 

  Let $B\subset T_0$ be a ball centred at $\tilde x$ such that $K=B\cap\partial C$ is contained in $\pmax C\,$. Let $H\subset S_0$ be the affine hyperplane through 
  $\tilde x$ which supports $C\,$. Since $\tilde x$ is exposed, there exists a half space $H_+$ containing 
  $\tilde x$ in its interior and with edge $H'$ parallel to $H\,$, such that $H_+\cap C\subset B\,$.

  Moreover, $H'\cap B\subset C$ is the base of a round (Lorentz) cone $W\subset T_0$ with apex $\tilde x$ and apex angle $0<\alpha<\pi\,$, and spans an affine hyperplane 
  parallel to the linear hyperplane $E_{1/2}(F)\subset T_0\,$. Its relative dual cone $W'$ with apex $0$ is a round cone with apex angle $\pi-\alpha\,$. Let $V$ be the 
  intersection of the relative interior of $W'$ with the unit sphere in $T_0\,$. Then, by elementary trigonometry, $p(V)$ is precisely the open ball around zero in $E_{1/2}(F)$ 
  with radius $\cos\tfrac\alpha2>0\,$. 

  We claim $V\subset\nu_0(K)\,$. In fact, if $y\in V\,$, then the affine hyperplane $H_y$ through $\tilde x$ with normal $y$ intersects $W$ in its apex $\tilde x\,$. 
  Therefore, $H_y\cap H'$ is disjoint from $B\cap H'=W\cap H'\,$. In other words, $H_y\cap C\subset H_+\,$. If $H_y=H\,$, then $y=\nu_0(\tilde x)\,$. Otherwise, $H_y$ 
  intersects $C$ in the interior of $H_+\cap C\,$, and the supporting hyperplane $H'_y$ of $C$ with inner normal $y$ also has $H'_y\cap C\subset H_+\,$. Thus, 
  $H_y'\cap C\subset B\cap\partial C=K\,$, and $y\in\nu_0(K)\,$. This proves the lemma. 
\end{pf}

\begin{rem}
 The article by Fedotov referred to in the proof of Lemma \ref{proj-open} is certainly known to experts in convexity, and R.~Schneider cites it in his monograph 
 \cite{schneider-convex}. 
\end{rem}

{\def\Elproofname{\protect{PROOF of Theorem \ref{c10-field}.}}
\begin{pf}
By Lemmata \ref{transverse} and \ref{proj-open}, $\mathbf e(\xi^{-1}(E))$ is a topological manifold of dimension $k=\dim E_{1/2}(F)=n_{d-j-1}-n_{d-j}-1\,$. Moreover, 
the injectivity of $p$ shows that $E_{1/2}(F)^\perp$ is transverse to $e_F(E)\in \mathbf e(\xi^{-1}(E))$ in the sense of Whitehead \cite[p.~155f.]{whitehead-transverse}. As 
explained there, $p:\mathbf e(\xi^{-1}(E))\to E_{1/2}(F)$ has a local Lipschitz inverse on some ball $B\subset E_{1/2}(F)\,$. Therefore, $\mathbf e(\xi^{-1}(E))$ is in 
fact a Lipschitz manifold. 

With respect to the Lipschitz chart $p^{-1}\,$, $\mathbf e(\xi^{-1}(E))$ has tangent spaces on a set of full $k$-dimensional Hausdorff measure in $p^{-1}(B)\,$. These 
tangent spaces have dimension $k\,$, and thus, by Lemma \ref{shilov-tangent}, in the $e_G(E)\in p^{-1}(B)\,$, the tangent space is $E_{1/2}(G)\,$. But the map 
$e_G(E)\mapsto E_{1/2}(G)$ is well defined and continuous with values the Grassmannian of $k$-planes in $\Span0E\,$. By the argument 
given at the beginning of this subsection, it follows that $p^{-1}$ is a regular $\mathcal C^1$ map. Since $p$ and the derivative of $p^{-1}$ depend continuously on $(E,F)\,$, it follows that $\mathcal P_j$ is indeed a closed $\mathcal C^{1,0}$ submanifold. 
\end{pf}}

\subsection{Construction of a $\mathcal C^{1,0}$ Groupoid and a Proper Homomorphism}

  We now proceed to the second step in the proof of our index theorem, as explained above. Since $\xi:\mathcal P_j\to\xi(\mathcal P_j)\subset P_{j-1}$ is a $\mathcal C^{1,0}$ manifold over $\xi(\mathcal P_j)\,$, we may consider the vector bundle $\vrho:\Sigma_{j-1}\to P_{j-1}\,$, and apply the construction of Subsection \ref{c10-grp-class} to 
  \[
  Y=\xi(\mathcal P_j)\,,\,M=\mathcal P_j\,,\,p=\xi\,,\,\nda\,E=\Sigma_{j-1}|\xi(\mathcal P_j)\ ,
  \]
  in the notation used there. Thus, 
  \begin{align*}
    \mathcal D_j=\xi^*\Sigma_{j-1}\times_{\mathcal P_j}(\mathcal P_j\times_{P_{j-1}}\mathcal P_j)=\Set1{(E,u,F_1,F_2)}{(E,F_i)\in\mathcal P_j\,,\,u\in E^\perp}
  \end{align*}
  is a $\mathcal C^{1,0}$ groupoid over $\mathcal D_j^{(0)}=\mathcal P_j\,$. Its Lie algebroid is $\xi^*\Sigma_{j-1}\oplus T\mathcal P_j\,$. Observe that since 
  \[
  F^\perp=E^\perp\oplus E_{1/2}(F)\oplus\reals\cdot e_F\mathfa(E,F)\in\mathcal P_j\ ,
  \]
  there is an isomorphism
  \[
  A(\mathcal D_j)\times\reals\cong\eta^*\Sigma_j:(E,F,u\oplus v,r)\mapsto(E,F,u+v+r\cdot e_F)
  \]
  of topological vector bundles over $\mathcal P_j\,$, in particular, of topological group\-oids.

  We now consider the groupoid $\mathbb W\mathcal D_j$ and define 
  \[
  \vphi:\mathbb W\mathcal D_j\to\mathcal W_\Omega|(U_{j+1}\setminus U_{j-1})
  \]
  by
  \[
  \vphi(\tau)=
  \begin{cases}
    \Parens0{E,r_1,u+r_2e_{F_2}-r_1e_{F_1}}&
    \left\{\begin{aligned}[c]
    \tau&=(E,u,F_1,F_2,r_1,r_2-r_1)\\
    &\in\mathcal D_j\times(\reals\rtimes\reals)|\reals_{\sge0}\ ;
    \end{aligned}\right.\\ & \\
    (F,0,u+v+re_F)&
    \left\{\begin{aligned}[c]
      \tau&=(E,F,u\oplus v,\infty,r)\\
      &\in A(\mathcal D_j)\times\infty\times\reals\ .
    \end{aligned}\right.
  \end{cases}
  \]

\begin{prop}\label{strmor-tang-wh}
  The map $\vphi$ is a proper strict morphism.
\end{prop}

\begin{pf}
  Recall that $\mathcal W_{\Omega}\subset\clos\Omega\times X$ is a closed embedding, so we may check the continuity of $\vphi$ component-wise. Equally,
  $\mathbb W\mathcal D_j\subset\mathbb T\mathcal D_j\times\mathcal W_{\reals_{\sge0}}$ is a closed embedding. Now,
  \[
  \mathcal P_j\to P_{j-1}\times X:(E,F,v)\mapsto(E,e_F)
  \]
  is a $\mathcal C^{1,0}$ chart. The corresponding $\mathcal C^{1,0}$ map on $\mathcal D_j$ is given by
  \[
  f:\mathcal D_j\to\mathcal P_j\times X:(E,u,F,G)\mapsto\Parens0{E,F,u+e_G-e_F}\ .
  \]

  Hence, let $(E,F,u\oplus v)\in \xi^*\Sigma_{j-1}\oplus T\mathcal P_j=A(\mathcal D_j)\,$. Observe
  \[
  r(E,F,u\oplus v)=(E,F)\equiv(E,0,F,F)\in\mathcal D_j\ .
  \]
  Since $E_{1/2}(F)=T_{e_F}\mathbf e(\xi^{-1}(E))\,$, there exist $F_\eps\in P_j\,$, $F_\eps\subset E\,$, \scth
  \[
  F=\limk_{\eps\to0+}F_\eps\nd v=\limk_{\eps\to0+}\frac{e_{F_\eps}-e_F}\eps\ .
  \]
  Moreover, $\gamma(\eps)=(E,\eps\cdot u,F,F_\eps)$ is a $\mathcal C^1$ curve in $r^{-1}(E,F,v)$ representing the tangent vector $(E,F,u+v)\,$. We find
  \[
  T_{(E,F)}f(u+v)=(f\circ\gamma)'(0)=(E,F,u+v)\ .
  \]
  By the definition of the topology on the tangent groupoid, the second component of $\vphi$ is continuous.

  As to the continuity of the first component, we need to see that
  \[
  -F^*=\limk_{\eps\to0+}\frac{e_{F_{1\eps}}}\eps-E_\eps^*\text{ in }\clos\Omega
  \]
  if $(E,F,u\oplus v)=\limk_{\eps\to0+}(E_\eps,u_\eps,F_{1\eps},F_{2\eps},\eps)\text{ in }\mathbb T\mathcal D_j\,$.  In fact, we have already seen that $F=\limk_{\eps\to0+}F_{1\eps}\,$, and $E=\limk_{\eps\to0+}E_\eps\,$. By continuity of polarity 
  \cite[Corollary 7.2.12]{beer-topologies-book}, $E^*=\limk_{\eps\to0+}E_\eps^*\,$. Then
  \begin{align*}
    \limk_{\eps\to0+}\frac{e_{F_{1\eps}}}\eps-E_\eps^*&=\limk_{\eps\to0+}\frac{e_F}\eps-E^*\\
    &=\limk_{\lambda\to\infty}\lambda\cdot e_F-E^*=\clos{\reals_{\sge0}\cdot e_F-E^*}=-F^*\ ,
  \end{align*}
  because $e_F$ generates the relative dual face of $F$ in $E^{\circledast}\,$, and by \cite[Lemma 2.2.4]{alldridge-johansen-wh1}. Therefore, $\vphi$ is continuous, and it is
  trivial to check that it indeed is a homomorphism.

  To see that $\vphi$ is proper, note first that $\mathcal W_\Omega|Y_j$ is closed in $\mathcal W_\Omega|(U_{j+1}\setminus U_{j-1})\,$. For any 
  $K\subset\mathcal W_\Omega|Y_j\,$, we have $\vphi^{-1}(K)\subset A(\mathcal D_j)\times\reals\,$. But the restriction of $\vphi$ to this set is the composition 
  of the projection $\eta^*\Sigma_j\to\Sigma_j\,$, which is proper, with the closed embedding $\Sigma_j\subset\mathcal W_\Omega|Y_j\,$. Hence, 
  $\vphi^{-1}(K)$ is compact if $K$ is. 

  Similarly, if $K\subset\mathcal W_\Omega|(U_{j+1}\setminus U_{j-1})$ is compact and completely contained in $\mathcal W_\Omega|Y_{j-1}\,$, then there exists
  $1\sle R<\infty$ \scth
  \[
  \max\Parens0{\Norm0{\lambda(r(\omega))},\Norm0{\lambda(s(\omega))}}\sle R\mathfa\omega\in K\ ,
  \]
  where we recall from \cite[Theorem 9]{alldridge-johansen-wh1} that $\lambda:\bar\Omega\to X$ is defined by $\lambda(x-F^*)=x$ and its restriction to $Y_{j-1}$ is 
  continuous. Furthermore,
  \[
  L=\Set1{u\in X}{\exists\,E\in P_{j-1}\,,\,v_1,v_2\in E^\circledast\,:\,(E,v_1,u+v_2-v_1)\in K}
  \]
  is compact. Hence,
  \[
  \vphi^{-1}(K)\subset P_{j-1}\times L\times P_j\times P_j\times[0,R]\times[-R,R]
  \]
  is compact in $\mathbb W\mathcal D_j\,$.

  It remains to consider a sequence
  \[
  \omega_k=\Parens1{E_k,r_k^1\cdot e_{F_k^1},u_k+r_k^2e_{F_k^2}-r_k^1e_{F_k^1}}
  \]
  converging to $(F,0,u+w)\in\Sigma_j\,,\,u\perp F\,,\,w\in F^\circledast\,$, and to exhibit a subsequence of $(E_k,u_k,F_k^1,F_2^2,r_k^1,r_k^2-r_k^1)$ converging 
  to $(E,F,u\oplus v,\infty,r)$ for some $E\in P_{j-1}\,$, $F\subset E\,$, $u\in E_{1/2}(F)\,$, $r\in\reals\,$, \scth $w=v+r\cdot e_F\,$. 

  In fact, by compactness of
  $\mathcal P_j\times_{P_{j-1}}\mathcal P_j\,$, by passing to a subsequence, we may assume that $(E_k,F_k^1,F_k^2)\to(E,F',F'')\,$. Moreover, since
  \[
  \limk_kr_k^1\cdot e_{F_k^1}-E_k^*\to-F^*\nd\dim F<\dim E=\dim E_k\ ,
  \]
  the sequence $r_k^1\cdot e_{F_k^1}$ is unbounded by \cite[Lemma 13 (iii)]{alldridge-johansen-wh1}, so $r_k^1\to\infty\,$. Because $E_k\to E$ 
  and we have $\dim E_k=\dim E\,$, we obtain $E_k^\perp\to E^\perp\,$. Thus,
  \[
  u_k+r_k^2\cdot e_{F_k^2}-r_k^1\cdot e_{F_k^1}\to u+w\mathtxt{implies}u_k\to u\ ,
  \]
  whence in turn $r_k^2\cdot e_{F_k^2}-r_k^1\cdot e_{F_k^1}\to w\,$. Because
  \[
  \Norm0{r_k^2\cdot e_{F_k^2}-r_k^1\cdot e_{F_k^1}}^2=\Parens0{r_k^2-r_k^1\cdot\Rscp0{e_{F_k^2}}{e_{F_k^1}}}^2+\Parens0{r_k^1\cdot(1-\Rscp0{e_{F_k^2}}{e_{F_k^1}})}^2\ ,
  \]
  passing to a subsequence, we may assume $r_k^2-r_k^1\to r\in\reals\,$, and $e_{F_k^2}-e_{F_k^1}\to0\,$, so $F'=F''\,$. Now,
  \[
  -F^*=\limk_kr_k^1\cdot e_{F_k^1}-E_k^*=\limk_{\lambda\to\infty}\lambda\cdot e_{F'}-E^*\ .
  \]
  This implies that $e_{F'}$ lies in the relative interior of $E^\circledast\cap F^\perp=\reals_{\sge0}\cdot e_F\,$, and hence $e_F=e_{F'}\,$, which finally gives
  $F=F'\,$. Since $(r_k^2-r_k^1)\cdot e_{F_k^2}\to r\cdot e_F\,$, the limit $v=\limk_kr_k^1\cdot(e_{F_k^2}-e_{F_k^1})$ exists. Necessarily, $v\in E_{1/2}(F)\,$, by Lemma 
  \ref{shilov-tangent}. We conclude $w=v+r\cdot e_F\,$. Thus, we have established the required relation
  \[
  (E_k,u_k,F_k^1,F_k^2,r_k^1,r_k^2-r_k^1)\to(E,F,u\oplus v,\infty,r)\mathtxt{in}\mathbb W\mathcal D_j\ ,
  \]
  and thereby, that $\vphi$ is proper.
\end{pf}

\subsection{Proof of the Main Theorem}

  As a corollary of the construction of the proper strict morphism $\vphi\,$, we obtain the topological Wiener--Hopf index formula on the level of operator $KK$ theory. To that end, let
  \[
  \vphi_0:A(\mathcal D_j)\times\reals=\eta^*\Sigma_j\to\Sigma_j\subset\mathcal W_\Omega|Y_j
  \]
  and
  \[
  \vphi_1:\mathcal D_j\times(\reals\rtimes\reals)|\reals_{\sge0}\to\mathcal W_\Omega|Y_{j-1}
  \]
  be the corresponding restrictions of $\vphi\,$. 

\begin{thm}\label{ind-opkk}
  We have the following expression for $\partial_j\,$:
  \[
  \partial_j\otimes KK(\vphi_1)=KK(\vphi_0)\otimes y\otimes\tau_j\mathtxt{in}KK^1\Parens1{\Cred0{\mathcal W_\Omega|Y_j},\Cred0{\mathcal D_j}}\ .
  \]
  Here, $\tau_j$ is the Connes--Skandalis map associated to the tangent groupoid $\mathbb T\mathcal D_j\,$, and the element 
  $y\in KK^1(S,\cplxs)$ is associated to the classical Wiener--Hopf extension.
\end{thm}

\begin{pf}
  Consider the strict morphism $\vphi:\mathbb W\mathcal D_j\to\mathcal W_\Omega|(U_{j+1}\setminus U_{j-1})$ from Proposition \ref{strmor-tang-wh}.
  Applying Corollary \ref{conn-hom-nat-kk}, we obtain
  \[
  \partial_j\otimes KK(\vphi_1)=KK(\vphi_0)\otimes\partial\ ,
  \]
  where $\partial$ represents the extension for $\mathbb W\mathcal D_j$ from Corollary \ref{wg-ext}. Now, the assertion follows from Theorem \ref{ytensortau}.
\end{pf}

  Consider the embeddings 
  \begin{gather*}
    i_{\mathcal P_j}:\mathcal P_j\to P_{j-1}\times X^2:(E,F)\mapsto (E,e_F,0)\ ,\\
    i_{\Sigma_{j-1}}:\Sigma_j\to P_{j-1}\times X^2:(E,u)\mapsto(E,u,0)\ ,
  \end{gather*}
  and the associated topological family index $\vphi_{T\mathcal P_j\oplus\xi^*\Sigma_{j-1}}^{\vphantom{-1}}\otimes\vphi_{\Sigma_{j-1}|\xi(\mathcal P_j)}^{-1}$ 
  (cf.~Section \ref{topexpr-connesskand}). 
  
  Recall from the introduction or \cite[Proposition 16]{alldridge-johansen-wh1} that the inclusion of $\Sigma_j$ in $\mathcal W_\Omega|Y_j$ is a Morita equivalence, and 
  similarly for $\Sigma_{j-1}\,$. By this token, the element 
  $\partial_j\in KK^1(\Cred0{\mathcal W_\Omega|Y_j},\Cred0{\mathcal W_\Omega|Y_{j-1}})$ may be pulled back to an element of 
  $KK^1(\Cred0{\Sigma_j},\Cred0{\Sigma_{j-1}})$ which we denote by the same letter. Applying Theorem \ref{cs-top-expr}, we obtain the following corollary.

\begin{cor}
  Let $\pmb\eta$ be the $*$-morphism induced by the projection $\eta^*\Sigma_j\to\Sigma_j$ (which is proper), and similarly let $\pmb\zeta$ be induced by the closed 
  embedding $\zeta:\Sigma_{j-1}|\xi(\mathcal P_j)\to\Sigma_{j-1}\,$. Then
  \[
  \pmb\zeta_*\partial_j=\pmb\eta^*\Bracks0{y\otimes\vphi_{T\mathcal P_j\oplus\xi^*\Sigma_{j-1}}^{\vphantom{-1}}\otimes\vphi_{\Sigma_{j-1}|\xi(\mathcal P_j)}^{-1}}
  \]
  in $KK^1\Parens0{\Cred0{\Sigma_j},\Cred0{\Sigma_j|\xi(\mathcal P_j)}}\,$. If $\xi:\mathcal P_j\to P_{j-1}$ is surjective, i.e.~every $n_{d-j+1}$-dimensional face of the 
  cone $\Omega^*$ contains an $n_{d-j}$-dimensional face, then $\pmb\zeta$ is the identity.
\end{cor}

\begin{pf}
  The assertion follows from Theorems \ref{ind-opkk} and \ref{cs-top-expr} by noting that $\eta^*=KK(\vphi_0)\,$, and that $\vphi_1$ drops to $\zeta$ through $\pi\,$.
\end{pf}

\end{document}